\documentclass[12pt]{article}

\usepackage{amsmath,mathtools,amssymb,amsthm,eucal,tikz,tikz-cd}
\usepackage[margin=1in]{geometry}
\usepackage[all,cmtip]{xy}

\usepackage{graphicx}
\usepackage[colorlinks=true,allcolors=blue,bookmarksopen,bookmarksdepth=3,linktocpage=true]{hyperref}
\usepackage[backrefs]{amsrefs}

\newcommand{\R}{\mathbb{R}}

\newcommand{\Z}{\mathbb{Z}}
\newcommand{\N}{\mathbb{N}}
\newcommand{\F}{\mathbb{F}}

\newcommand{\cfuk}{\mathcal{F}uk}

\DeclareMathOperator{\Hom}{Hom}
\DeclareMathOperator{\id}{id}

\DeclareMathOperator{\End}{End}

\DeclareMathOperator{\Mod}{Mod}
\DeclareMathOperator{\Star}{Star}

\newcommand{\Sh}{\operatorname{Sh}^{w}}

\numberwithin{equation}{section}

\theoremstyle{plain}
\newtheorem{thm}[equation]{Theorem}
\newtheorem*{thm*}{Theorem}
\newtheorem{cor}[equation]{Corollary}
\newtheorem{lem}[equation]{Lemma}
\newtheorem{prop}[equation]{Proposition}

\theoremstyle{definition}
\newtheorem{defn}[equation]{Definition}

\theoremstyle{remark}
\newtheorem{rem}[equation]{Remark}

\newtheorem{assumption}[equation]{Assumption}

\begin{document}

\title{Partially wrapped contact Fukaya categories\\ of 1-jet spaces}

\author{Hanming Liu\footnote{Department of Mathematics, University of Oregon. Email: \textit{hanming@uoregon.edu}}}

\date{November 1, 2026}

\maketitle

\begin{abstract}
	\noindent We introduce the partially wrapped contact Fukaya category of a \(1\)-jet space equipped with a stop, which is an \(A_\infty\)-category of Legendrians equipped with representations of their Chekanov-Eliashberg algebras. This category is invariant under contact isotopies of the stop. We prove a wrapping exact triangle and a generation theorem. When the stop is bounded subanalytic singular isotropic, we prove a Morita equivalence between the partially wrapped contact Fukaya category and a microlocal sheaf category.
\end{abstract}

\newpage
\tableofcontents
\newpage

\section{Introduction}

In 1993, Fukaya \cite{Fukaya} introduced an \(A_\infty\)-category structure on Lagrangian submanifolds in a symplectic manifold \(W\), later known as the Fukaya category of \(W\). Fukaya categories are constructed using the theory of holomorphic curves, first developed by Gromov \cite{Gromov} and Floer \cites{Floer1,Floer2,Floer3,Floer4}. As Floer theory counts solutions to partial differential equations, Fukaya categories are analytical in nature. In 1994, Kontsevich \cite{Kontsevich} introduced the Homological Mirror Symmetry conjecture, which relates the Fukaya category of a symplectic manifold with the derived category of sheaves of an algebraic variety. This conjecture has shaped symplectic topology and algebraic geometry ever since.

Microlocal sheaf theory, developed by (among others) Kashiwara-Schapira \cites{KS}, has also been a fruitful tool for studying symplectic topology, for example \cites{Tamarkin,GKS}. Computations in microlocal sheaf theory are combinatorial and algebraic in nature.

Fukaya categories and microlocal sheaf categories turn out to be closely related. Nadler-Zaslow \cite{NZ} and Nadler \cite{Nadler} proved that the infinitesimal Fukaya category of a cotangent bundle \(T^*M\) is equivalent to the category of sheaves with perfect stalks. Later, Ganatra-Pardon-Shende \cite{GPS3} proved an analogous equivalence between partially wrapped Fukaya categories and the category of compact sheaves with appropriate microsupport. These equivalences unite certain aspects of analytical and combinatorial approaches to symplectic topology.

Contact topology is the odd-dimensional sibling of symplectic topology. Both holomorphic curve techniques (for example \cites{EHS,EGH,Che,Ng,EES3,Pardon}) and sheaf theoretic techniques (for example \cites{Shende,CG1,CG2,GSW,Li,Kilgore}) have proved to be fruitful in contact topology.

In analogy to Fukaya categories, a \emph{contact Fukaya category} should heuristically be a category with objects being Legendrian submanifolds in a contact manifold. In this paper, we introduce a version of contact Fukaya categories. For a \(1\)-jet space \(J^1M\) and a stop \(\mathfrak{f}\subset J^1M\) (which can be any closed subset), we construct the \emph{partially wrapped contact Fukaya category} \(\mathcal{V}(J^1M,\mathfrak{f})\) (Definition \ref{definition of partially wrapped contact Fukaya categories}), which is an \(A_\infty\)-category of \emph{Legendrian \(R\)-branes} \(\Lambda\subset J^1M\setminus \mathfrak{f}\); A Legendrian \(R\)-brane is a Legendrian equipped with some extra data, including a representation of its Chekanov-Eliashberg algebra (see Definition \ref{definition R-brane}). Morphisms in \(\mathcal{V}(J^1M,\mathfrak{f})\) are cochain complexes whose cohomology groups are given by the \emph{partially wrapped bilinearized contact cohomology}. Writing \(\mathcal{V}(\Lambda_1,\Lambda_2)\) for \(\Hom_{\mathcal{V}(J^1M,\mathfrak{f})}(\Lambda_1,\Lambda_2)\), this partially wrapped bilinearized contact cohomology is the colimit \[H^*\mathcal{V}(\Lambda_1,\Lambda_2):=\varinjlim_i LCH(\Lambda_1,\Lambda_2^t),\] where \(\Lambda_2^t\) is a positive cofinal wrapping of \(\Lambda_2\) in \(J^1M\setminus\mathfrak{f}\) (see \cite{GPS2}*{Section 2.1}), and \(LCH\) is the \emph{bilinearized contact cohomology} (Section \ref{section continuation elements}). The construction of \(\mathcal{V}(J^1M,\mathfrak{f})\) is inspired by earlier work on the augmentation category \cites{BC,NRSSZ,CDGG,CNS,Liu25} and the partially wrapped Fukaya category \cites{Sylvan,GPS1,GPS2}.

By a \emph{contact isotopy}, we mean a family of contactomorphisms starting from the identity.

\begin{thm*}
   (Proposition \ref{proposition well-defined}, Theorem \ref{isotopy invariance}) Up to quasi-equivalence, \(\mathcal{V}(J^1M, \mathfrak{f})\) is well-defined, and invariant under contact isotopies of \(\mathfrak{f}\).
\end{thm*}

When the stop is Legendrian, a wrapping exact triangle quantifies the effect of wrapping a Legendrian through the stop, which is the main technical ingredient in the computation of partially wrapped contact Fukaya categories:

\begin{thm*}
   (Theorem \ref{wrapping exact triangle}) Let \(p\in \mathfrak{f}\) be a point around which locally \(\mathfrak{f}\) is a Legendrian submanifold of \(V\). Let \(\Lambda\) be a Legendrian \(R\)-brane, and let \(\Lambda^w\in V\) be a Legendrian \(R\)-brane such that there is a positive isotopy from \(\Lambda\) to \(\Lambda^w\in V\) in \(V\setminus(\mathfrak{f}\setminus p)\) which goes through \(p\) once. Let \(S\) be a Legendrian linking sphere of \(\mathfrak{f}\) at the point \(p\). Then, in the category \(\mathcal{V}(J^1M,\mathfrak{f})\), there is an exact triangle
	\[
	\Lambda\rightarrow\Lambda^w\rightarrow S.
	\]
\end{thm*}

Another main ingredient in the computation of partially wrapped contact Fukaya categories is a characterization of \(\mathcal{V}(J^1M,\emptyset)\):

\begin{thm*}
   (Theorem \ref{compact wcfuk with empty stop}) \(\mathcal{V}(J^1M,\emptyset)\) is equivalent to the compact Fukaya category of \(T^*M\), which is also equivalent to the derived category of finite dimensional local systems on \(M\), and is generated by Legendrian \(R\)-branes whose underlying Legendrian is a section of \(J^1M\).
\end{thm*}

There is an analogously defined category of sheaves constructed by Kuo \cite{Kuo}, the partially wrapped finite rank microlocal sheaf category, which we denote \(\Sh(J^1M,\mathfrak{f})\). The partially wrapped finite rank microlocal sheaf category is built out of sheaves of \(M\times\R\) microsupported on a compact Legendrian in \(J^1M\setminus\mathfrak{f}\) (which we call \emph{Legendrian \(S\)-branes}), and it satisfies structural theorems that are analogous to the two theorems above (Theorem \ref{sheaf wrapping exact triangle}, Theorem \ref{sheaf with empty stop}). This sheaf category turns out to be equivalent to the partially wrapped contact Fukaya category:

\begin{thm*}\label{intro main theorem}
	(Theorem \ref{main theorem}) Let \(M\) be a real analytic manifold. Let \(\mathfrak{f}\) be a subanalytic singular isotropic stop bounded in the vertical direction. If \(M\) is 4- or 5-dimensional, we additionally assume Assumption \ref{assumption}. Then, there is a Morita equivalence \[\mathcal{V}(J^1M,\mathfrak{f})\cong \Sh(J^1M,\mathfrak{f}).\]
\end{thm*}

The proof strategy of Theorem \ref{main theorem} roughly goes as follows. We first enlarge both categories to allow for some infinite-dimensional local systems. Then we enlarge \(\mathfrak{f}\) to be the conormal of a subanalytic Whitney stratification of \(M\times\R\). Then both (enlarged) categories have a set of generators for which we can write down explicit wrapping sequences. If the stratification is chosen carefully enough, then morphisms between generators are concentrated in degree \(0\), and compositions are easy to compute. It turns out that generators, morphisms between generators, and their compositions in both enlarged categories match up, thereby giving a Morita equivalence between the enlarged categories. We then use the \(A_\infty\)-Yoneda lemma together with the stop removal theorem (Theorem \ref{stop removal}) to prove that the original categories are equivalent. Essentially, we find generators and relations for both sides, and prove that they are the same.

Theorem \ref{main theorem} is related to the augmentation-sheaf correspondence first proposed in \cite{STZ}, proven in \cite{NRSSZ} for Legendrian links in \(\R^3\) (later developments include \cites{CNS,ABS}): \cite{NRSSZ} proves a quasi-equivalence between the augmentation category \(\mathcal{A}ug_+(\Lambda)\) of a Legendrian link \(\Lambda\subset \R^3\) and the rank 1 microlocal sheaf category \(\operatorname{Sh}^1(\Lambda)\). By construction, \(\mathcal{A}ug_+(\Lambda)\) is a full subcategory of \(\mathcal{V}(\R^3,\Lambda)\), and \(\operatorname{Sh}^1(\Lambda)\) is a full subcategory of \(\Sh(\R^3,\Lambda)\). The augmentation-sheaf correspondence does not generalize to all \(1\)-jet spaces: it fails in \(J^1S^2\) \cites{Gao,GL}. As such, we view Theorem \ref{main theorem} as a broader framework in which the augmentation-sheaf correspondence of \(\R^3\) lives, and is the correct generalization of the Floer-sheaf correspondence from \(\R^3\) to arbitrary \(1\)-jet spaces.

The paper is organized as follows. In Section \ref{section preliminaries}, we recall the basics of the Chekanov-Eliashberg algebra, and recall (and reformulate) some algebraic preliminaries on augmented DGAs and \(A_\infty\)-categories to be used later. In Section \ref{section partially wrapped contact Fukaya categories}, we construct the partially wrapped contact Fukaya categories, and prove some structural results including the wrapping exact triangle. In Section \ref{section sheaf categories}, we record some structural results about the partially wrapped microlocal sheaf categories. In Section \ref{section equivalence}, we prove the equivalence between the partially wrapped contact Fukaya category and the partially wrapped microlocal sheaf category.

\section*{Acknowledgements}

I thank my advisor Robert Lipshitz for his guidance and support. I thank Honghao Gao, Andrew Hanlon, Lenhard Ng, John Pardon, Harahm Park, Dan Rutherford, Vivek Shende for helpful discussions and correspondences. I thank Tobias Ekholm, Zijun Li, Jiajie Ma for helpful discussions regarding orientations of moduli spaces. I thank Wenyuan Li for his numerous insights which have been extremely helpful during times when this project was stuck. I especially thank Frank Zheng, who even though generously declined to be named as a coauthor, has discussed this project with me extensively throughout its entirety, and contributed significantly to several proofs. The write-up of this paper was finished during my stay at the Institute for Pure and Applied Mathematics in fall 2026, for which I thank their hospitality. I was supported by NSF Grant DMS-2505715 and a Simons dissertation fellowship.

\section{Preliminaries}\label{section preliminaries}

In this section, we introduce the Chekanov-Eliashberg differential graded algebra, and the \(A_\infty\)-structure coming from dualizing an augmented differential graded algebra. Nothing in this section contains essentially new ideas, but our treatment differs somewhat from existing literature \cites{EES3,NRSSZ,CNS}.

In section \ref{section cedga}, we discuss the Chekanov-Eliashberg differential graded algebra, and introduce representations of Legendrians. Sections \ref{section link gradings and split DGAs}, \ref{section augmented DGAs}, \ref{algebraic functoriality} are devoted to algebraic formalism behind the proofs of \(A_\infty\)-relations and invariance of partially wrapped contact Fukaya categories, which appear in Section \ref{section partially wrapped contact Fukaya categories}. Since we state a formula for the \(A_\infty\)-operations in Section \ref{section unwrapped}, the reader might prefer to skip Sections \ref{section link gradings and split DGAs}, \ref{section augmented DGAs}, \ref{algebraic functoriality} on the first pass.

\subsection{Chekanov-Eliashberg algebras and representations}\label{section cedga}

We review some background about the Chekanov-Eliashberg differential graded algebra, and representations. We refer the reader to \cite{EES3} for details, and \cite{EENS} for another exposition on the treatment for disconnected Legendrian submanifolds. We use a slightly unconventional coefficient ring, and as a result, our definition of representations is uncommon as well (different from, say, \cite{CNS}). Our choice of coefficient ring is to make it easy for us to talk about representations of disconnected Legendrians where the dimension of the representation varies across components.

Let \((P,d\lambda)\) be a \(2n\)-dimensional exact symplectic manifold with finite geometry at infinity. See \cite{EES3}*{Definition 2.1} for the notion of finite geometry at infinity. Let \(V=P\times\R\). Equip \(V\) with the standard contact form \(\alpha=dz-\lambda\), where \(z\) is the \(\R\) component (\(V\) is sometimes called the \textit{contactization} of \(P\)). Let \(\Pi:V\to M\) denote the Lagrangian projection.

With the contact form \(\alpha\), the Reeb vector field \(R_\alpha\) satisfies \(R_\alpha=\partial_z\). So, given a Legendrian submanifold \(\Lambda\subset V\), the Reeb chords of \(\Lambda\) correspond to self intersections of \(\Pi(\Lambda)\). We can perturb \(\Lambda\) so that there are finitely many Reeb chords, and the corresponding self intersections of \(\Pi(\Lambda)\) are transverse double points. If the Legendrian \(\Lambda\) satisfies this condition, we call it \emph{chord generic}. We denote the set of Reeb chords by \(\mathcal{R}\).

\begin{defn}
Let \(X\subset Y\) be a submanifold. \(X\) is said to be \emph{relatively pin} with respect to a cohomology class \(b\in H^2(Y;\Z/2)\), if the restriction of \(b\) to \(X\) is the second Stiefel-Whitney class \(w_2(X)\).
\end{defn}

Let \(\Lambda\subset V\) be a chord generic compact relatively pin Legendrian submanifold with vanishing Maslov class. We fix a point \(p_j\in\Lambda_j\) on each component \(\Lambda_j\) of \(\Lambda\). We pick one of the basepoints, say \(p_1\) for convenience, which we call the \emph{special basepoint}. For each Reeb chord endpoint in \(\Lambda_j\) we choose an \emph{endpoint path} connecting the endpoint to \(p_j\). Furthermore, for \(j=1,...,r\), we choose paths \(\gamma_{1j}\) in \(P\) connecting \(\Pi(p_1)\) to \(\Pi(p_j)\) and symplectic trivialisations of \(\gamma_{1j}^*TP\) in which the tangent space \(\Pi(T_{p_1}\Lambda)\) corresponds to the tangent space \(\Pi(T_{p_j})\); for \(j=1\), \(\gamma_{11}\) is the trivial path. For any \(1\leq i,j\leq r\), we can then define \(\gamma_{ij}\) to be the path \(\gamma_{1i}\cup \gamma_{1j}\) joining \(\Pi(p_i)\) to \(\Pi(p_j)\), and \(\gamma_{ij}^*TM\) inherits a symplectic trivialisation from the trivialisations for \(\gamma_{1i}\) and \(\gamma_{1j}\).

Let \(R:=\bigoplus_{i}\Z[\pi_1(\Lambda_i,p_i)]\) be the direct sum of the group rings of the fundamental groups of all connected components of \(\Lambda\).

\begin{defn}\label{defn of CEDGA}
	The \textit{Chekanov Eliashberg differential graded algebra}, \((CE(V,\Lambda),\partial)\) is defined as follows:
	\begin{enumerate}
		\item \textbf{Algebra}: \(CE(V,\Lambda)\) is the free non-commutative algebra over the ring \(R\), generated by elements of \(\mathcal{R}\), where the coefficients do not commute with the generators. In other words, \(CE(V,\Lambda)\) is the non-commutative unital ring generated by elements of \(\mathcal{R}\) and \(\pi_1(\Lambda_i)\), with relations of the form \(hg=h\cdot_{\pi_1(\Lambda_i)}g\) for all \(h,g\in \pi_1(\Lambda_i)\), and the relation that the identity of \(R\) is identified with \(1\).

		\item \textbf{Grading}: Elements of \(R\) have degree 0. If \(q\) is a Reeb chord we define the grading by associating a path of Lagrangian subspaces to \(q\). If the endpoints of \(q\) lie on the same component \(\Lambda_j\), consider the path of tangent planes along the endpoint path from the final point of \(q\) to \(p_j\) followed by the reverse endpoint path from \(p_j\) to the initial point of \(q\). We then close this path to a loop \(\hat{\gamma}_q\) (see \cite{EES1} for details). In the case where the endpoints of \(q\) lie on different components, we associate a loop of Lagrangian subspaces \(\hat{\gamma}_q\) to \(q\) in the same way except that we insert the path of Lagrangian subspaces induced by the trivialisation along the chosen path, \(\gamma_{ij}\), connecting the components of the endpoints in order to connect the two paths from Reeb chord endpoints to chosen points. The grading of \(q\) is then \[|q|=\mu(\hat{\gamma}_q)-1,\] where \(\mu\) denotes the Maslov index. Note that the grading is independent of choices of endpoint paths in the case where \(\Lambda\) is connected, and the grading depends on the choices of paths in the case where \(\Lambda\) is disconnected.
		
		In the case where \(P\times\R=J^1M\) is a \(1\)-jet space, the choice of these paths is equivalent to choosing a Maslov potential. A Maslov potential \(\mu\) on \(\Lambda\) assigns an integer to points on \(\Lambda\), in a way that is locally constant, except for at cusps of the front projection, where the integer increases by one when going up through a cusp. Denote the Maslov potential of the beginning point (respectively end point) of a Reeb chord \(q\) by \(\mu_0(q)\) (respectively \(\mu_1(q)\)). In the front projection, away from cusps, a Legendrian is a multi-function \(f:M\to\R\), see \cite{Ekh}*{Section 2.2.1} for detail. A Reeb chord \(q\) corresponds to a critical point of the difference of the different sheets of \(f\), and we denote the index of this critical point by \(\operatorname{ind}(q)\). The grading is then given by the formula
        \[
        |q|=\operatorname{ind}(q)+\mu_1(q)-\mu_0(q)-1.
        \]

		\item \textbf{Differential}: We define the differential on generators, and extend by the Leibniz rule. We set \(\partial (R)=0\). For a Reeb chord \(q_0\), \[\partial q_0=\sum_{{
			\begin{matrix}
				&\bar{A}=(A_0,\dots,A_k),\\
				&\sum_{j=1}^k|q_j|=|q_0|-1\\
			\end{matrix}}}
		(-1)^{(n-1)(|q_0|+1)}|\mathcal{M}_{\bar{A}}(q_0;q_1...q_k)|A_0q_1A_1q_2...A_{k-1}q_kA_k,\] where \(A_i\) is an element of \(\pi_1(\Lambda_j)\) for some \(j\), and \(q_i\) are Reeb chords, and \(|\mathcal{M}_{\bar{A}}(q_0;q_1,...,q_k)|\) is the signed count of the moduli space defined below.
	\end{enumerate}
\end{defn}

Now we define the moduli space \(\mathcal{M}_{\bar{A}}(q_0;q_1...q_k)\). Let \(J\) be an almost complex structure on the symplectisation \((V\times\R,d(e^t\alpha))\) of \(V\) (where \(\alpha\) is the contact form we equipped \(V\) with, and \(t\) is the \(\R\) coordinate) that is compatible with the symplectisation in the following sense: \(J\) is \(\R\)-invariant, \(J(\partial/\partial t)=R_\alpha\), and \(J\) maps \(\xi=\ker\alpha\) to itself. With respect to this almost complex structure, \(q_i \times \R\) is a holomorphic strip for any Reeb chord \(q_i\) of \(\Lambda\).

Let \(D_k^2=D^2\setminus \{p^+,p^-_1,...,p^-_k\}\) be a closed disk with \(k+1\) punctures on its boundary, labeled \(p^+,p^-_1,...,p^-_k\) in counterclockwise order around \(\partial D^2\).

\begin{figure}
   \centering
	\includegraphics{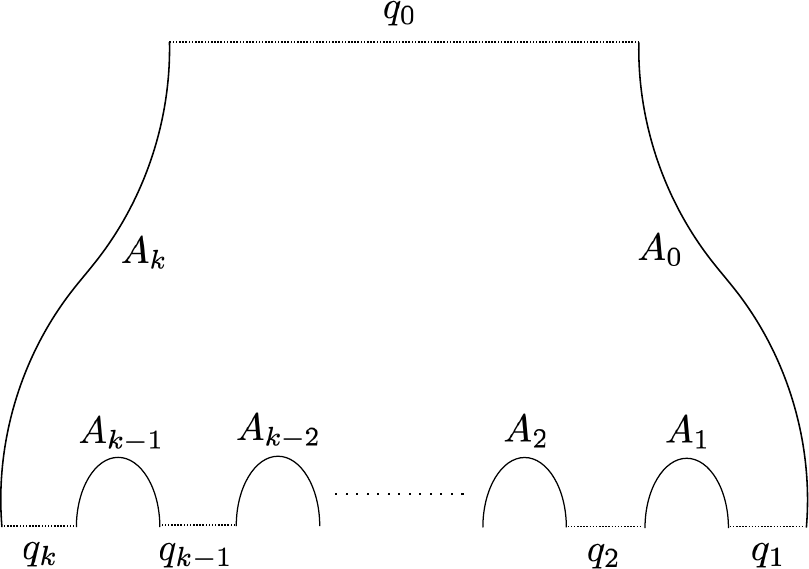}
	\caption{Holomorphic disks counted in the Chekanov-Eliashberg algebra}
	\label{figure DGA disk}
\end{figure}

For a \((k+1)\)-tuple of loops \(\bar{A}=(A_0,\dots,A_k)\) (more precisely, the loops are elements of \(\pi(\Lambda_j)\) for some \(j\)), we let \[\mathcal{M}_{\bar{A}}(q_0;q_1,...,q_k)\] denote the moduli space of \(J\)-holomorphic maps \[u:(D^2_k,\partial D^2_k)\to (V\times\R,\Lambda\times\R)\] up to domain reparametrisation and \(\R\)-translation, such that:
\begin{enumerate}
	\item near \(p^+\), \(u\) is asymptotic to a neighborhood of the Reeb strip \(q_0\times\R\) near \(t=\infty\);
	\item near \(p_l^-\) for \(1\leq l\leq k\), \(u\) is asymptotic to a neighborhood of \(q_l\times\R\) near \(t=-\infty\);
	\item when \((\Pi\circ u)|_{(p_j,p_{j+1})}\) is completed to a loop using the endpoint paths, it represents the class \(A_j\). Here \((p_j,p_{j+1})\) denotes the boundary interval in \(\partial D^2\) between \(p_j\) and \(p_{j+1}\).
\end{enumerate}
See Figure \ref{figure DGA disk} for a picture of the holomorphic disks in \(\mathcal{M}_{\bar{A}}(q_0;q_1,...,q_k)\).

\begin{rem}
In \cite{EES2}, the orientation of the moduli space \(\mathcal{M}_{\bar{A}}(q_0;q_1,...,q_k)\) depends on the choice of a spin structure on \(\Lambda\). By an argument similar to that of \cite{Abouzaid4}*{Lemma 6.1}, the spin structure can be relaxed to a relative pin structure.
\end{rem}

\begin{defn}
We will call such a \(J\)-holomorphic map a holomorphic disk (or simply disk) on \(\Lambda\) with a positive puncture at \(q_0\) and negative punctures at \(q_1,\ldots,q_k\).
\end{defn}

It is proved in \cite{EES3} that for generic admissible \(\Lambda\) and \(J\), \(\mathcal{M}_{\bar{A}}(q_0;q_1,...,q_k)\) is a transversely cut out manifold of dimension \(|q_0|-\sum_{i=1}^{k}|q_i|-1\). Roughly, admissibility here means that around the endpoints of Reeb chords, \(\Lambda\) is real analytic with respect to some standard coordinate given by \(J\), and generic admissible means that we further perturb \(\Lambda\) in the space of admissible Legendrian submanifolds (for the precise notion of generic admissible, see the discussion surrounding \cite{EES3}*{Section 2.3}). Furthermore, it admits a compactification as a manifold with boundary with corners in which the boundary consists of broken disks. Consequently, if the dimension equals \(0\) then it is compact.

Next we review the notion of stable tame isomorphism, first introduced by Chekanov in \cite{Che}. Given a semi-free DGA \(\mathcal{A}=R\langle q_1,\ldots,q_n\rangle\), an \emph{elementary automorphism} is an algebra automorphism of \(\mathcal{A}\) that sends a free generator \(q_i\) to \(xq_i y+u\) (here, \(x,y\in R\) are coefficients, and \(u\) is a linear combination of words not involving q), and sends all other generators to themselves. A \emph{tame automorphism} is a composition of elementary automorphisms. A \emph{tame isomorphism} between \((R\langle q_1,\ldots,q_n\rangle,\partial_1)\) and \((R\langle q'_1,\ldots,q'_n\rangle,\partial_2)\) is a chain map given by the composition of a tame automorphism and the map \(q_i\mapsto q'_i\) for all \(i\). A \emph{stabilisation} of \((R\langle q_1,\ldots,q_n\rangle,\partial)\) is \((R\langle q_1,\ldots,q_n,e_1,e_2\rangle,\partial)\) such that \(\partial(e_1)=e_2\), and \(\partial\) agrees with the original \(\partial\) on all other generators, and we also call the inclusion \((R\langle q_1,\ldots,q_n\rangle,\partial)\hookrightarrow (R\langle q_1,\ldots,q_n,e_1,e_2\rangle,\partial)\) a stablisation. A \emph{destabilisation} is the quotient map \((R\langle q_1,\ldots,q_n,e_1,e_2\rangle,\partial)\twoheadrightarrow (R\langle q_1,\ldots,q_n\rangle,\partial)\). A \emph{stable tame isomorphism} is a composition of tame isomorphisms, stabilisations, and destabilisations.

\begin{thm}[\cite{EES3}]\label{EES theorem}
	The map \(\partial:CE(V,\Lambda)\to CE(V,\Lambda)\) is a differential, that is \(\partial^2=0\). The stable tame isomorphism class of \(CE(V,\Lambda)\) is an invariant of \(\Lambda\) up to Legendrian isotopy.
\end{thm}

\begin{rem}
	In \cite{EES3}, they proved the version where the ring \(R\) is replaced by \(\Z[H_1(\Lambda)]\), and where homology classes commute with Reeb chords. The proofs carry over to our setting with no change.
\end{rem}

\begin{defn}
   Let \(\F\) be a field. Let \(\Lambda\subset V\) be a compact relatively pin Legendrian with vanishing Maslov class, equipped with data needed to define a grading on the Reeb chords. Write \(\Lambda_i\) for the connected components of \(\Lambda\). For each component \(\Lambda_i\), let \(V_i\) be a finite dimensional vector space over \(\F\). Write \(V:=\bigoplus_{i}V_i\). We view \(\End(\bigoplus_{i}V_i)=\bigoplus_{i,j}\Hom_{\F}(V_i,V_j)\) as a DGA concentrated in degree \(0\). A \emph{representation} \(\epsilon\) of \(\Lambda\) is a DGA map \[\epsilon:CE(V,\Lambda)\to \End(V),\] such that for each Reeb chord \(b_{ij}\) from \(\Lambda_i\) to \(\Lambda_j\), \(\epsilon(b_{ij})\in \Hom_{\F}(V_j,V_i)\), and that for each element \(\gamma_i\in\pi_1(\Lambda_i)\), \(\epsilon(\gamma_i)\in \End(V_i)\).

   In the language of Definition \ref{definition respect the link grading}, a representation of \(\Lambda\) is equivalently a representation of \(\mathcal{A}(\Lambda)\) that respects the link grading coming from decomposing \(\Lambda\) as its connected components.
\end{defn}

\subsection{Link gradings}\label{section link gradings and split DGAs}

We review the notion of link grading, which is simply a way to keep track of which components the endpoints of a Reeb chord lie on. We note that our definition of link gradings differ slightly from that of \cites{NRSSZ,Liu25}, see Remark \ref{remark link grading}.

\begin{defn}\label{semi-free dga}
	A \emph{semi-free} DGA is a DGA that as an algebra, is the free non-commutative algebra generated by a set \(\mathcal{R}\), over a non-commutative coefficient ring \(R\), such that elements of the coefficient ring \(R\) do not commute with the generators.
	
	In other words: let \(R\) be a non-commutative coefficient ring, with a generating set \(\mathcal{H}\), and some relations. A semi-free DGA is a DGA \(\mathcal{A}\) equipped with a set \(\mathcal{S}=\mathcal{R}\cup\mathcal{H}\) of generators, such that \(\mathcal{A}\) is the result of taking the free non-commutative unital algebra generated by the elements of \(\mathcal{S}\), and quotienting by relations on \(\mathcal{H}\) that are inherited from the relations that define \(R\). In addition, for any \(h\in\mathcal{H}\), we require that \(|h|=0\), and \(\partial h=0\).
\end{defn}

\begin{rem}
Examples of semi-free DGAs as defined above include Chekanov-Eliashberg DGAs for a pair \((V,\Lambda)\), setting \(\mathcal{R}\) to be the set of Reeb chords, and \(\mathcal{H}\) is a generating set of \(R\).
\end{rem}

\begin{defn}\label{definition of link grading}
	Let \(\mathcal{A}\) be a semi-free DGA with generating set \(\mathcal{S}=\mathcal{R}\cup\mathcal{H}\). Let \(P\) be a finite set. A \(P\)-\textit{link grading} on \(\mathcal{A}\) is a choice of a pair of maps \[c,r:\mathcal{S}\to P\] satisfying the following conditions:
	\begin{enumerate}
		\item for any \(x\in\mathcal{S}\), each term in \(\partial x\) is a word of the form \(x_1...x_k\) where \(c(x_{i})=r(x_{i+1})\) for \(i=1,...,k-1\) and \(c(x_k)=c(x),r(x_1)=r(x)\) (such a word is called \textit{composable} from \(c(x)\) to \(r(x)\));
		\item for any \(h\in\mathcal{H}\), we have \(c(h)=r(h)\).
	\end{enumerate}
\end{defn}

\begin{rem}\label{remark link grading}
   Note that our definition of link gradings here differs slightly from the definitions found in \cites{Liu25,NRSSZ}. The reason is that we will be using this definition along with our definition of the Chekanov-Eliashberg algebra, which has a different coefficient ring than \cites{Liu25,NRSSZ}. Our definition of link gradings here (along with the definition of link gradings in \cite{Liu25}) is also different from that of \cite{NRSSZ} in that the superscripts are of the opposite order to theirs, since in this paper we want the composition in the partially wrapped contact Fukaya category to obey the forward composition convention (see Section \ref{section augmented DGAs}).
\end{rem}

\begin{defn}\label{Legendrian link grading}
	If \(\Lambda=\coprod_P \Lambda^i\) is a disjoint union of compact relatively pin Legendrians indexed by a finite set \(P\), each equipped with data needed to define the Chekanov-Eliashberg algebra, then we give \(\mathcal{A}(\Lambda)\) a \(P\)-link grading as follows: for each Reeb chord \(q\), define \(c(q)\) (respectively \(r(q)\)) to be the index of the component containing the beginning point (respectively endpoint) of \(q\), and for \(\gamma\) a based loop in some connected component of \(\Lambda^i\), define \(r(\gamma)=c(\gamma)=i\).
\end{defn}

Let \(\mathcal{A}\) be equipped with a \(P\)-link grading.

\begin{defn}
   For \(i,j\in P\), we write \(\mathcal{S}^{ij}:=(c\times r)^{-1}(i,j)\), and likewise \(\mathcal{R}^{ij}\) and \(\mathcal{H}^{ij}\). We call elements of \(\mathcal{S}^{ii}\) \textit{diagonal} and elements of \(\mathcal{S}^{ij}\) for \(i\neq j\) \textit{off-diagonal}. Note that all elements of \(\mathcal{H}\) are diagonal.
\end{defn}

\begin{defn}\label{definition quotient DGA}
   Let \(P'\subset P\) be a subset. Define \(\mathcal{A}_{P'}:=\mathcal{A}/\left(\bigcup_{i\in P'\text{ or }j\in P'}S^{ij}\right)\). Write \(\mathcal{A}_{i}\) for \(\mathcal{A}_{\{i\}}\).
\end{defn}

\begin{defn}\label{definition respect the link grading}
   Let \(\mathcal{A},\mathcal{A}'\) be DGAs with generating sets \(\mathcal{S},\mathcal{S}'\), both equipped with a \(P\)-link grading. A map of DGAs \(f:\mathcal{A}\to\mathcal{B}\) is said to \emph{respect the link grading} if \(f(\mathcal{S}^{ij})\) consists of composable words in \(\mathcal{A}'\) from \(i\) to \(j\).

   A representation \(\epsilon:\mathcal{A}\to\End(V)\) is said to respect the link grading if there exists a decomposition \(V=\bigoplus_{i\in P}V_i\), so that: if we equip \(\End(V)\) (as a DGA concentrated in degree \(0\)) with the link grading such that each non-zero element \(x\in\Hom(V_j,V_i)\) has link gradings \(c(x)=i,r(x)=j\), then \(\epsilon\) is a DGA map that respects the link grading.
\end{defn}

\begin{defn}
   A representation \(\epsilon:\mathcal{A}\to \End(V)\) is called \emph{diagonal} if off-diagonal elements of \(\mathcal{S}\) get sent to \(0\). Equipping \(\mathcal{A}\) with a diagonal representation is equivalent to equipping each \(\mathcal{A}_i\) with a representation \(\epsilon_i\).
\end{defn}

\subsection{\(A_\infty\)-operations}\label{section augmented DGAs}

We construct the \(A_\infty\)-operations that will be used to define the partially wrapped contact Fukaya category, roughly following the outline in \cite{CNS}*{Section 2.1}.

We work with \(\Z/2m\)-graded \(A_\infty\)-categories over a field \(\F\). We use the \emph{forward composition convention}: let \(\mathcal{C}\) be an \(A_\infty\)-category, and let \(X,Y\) be objects. Let \(\mathcal{C}(X,Y)\) denote the hom space. We write composition as \[m_k:\mathcal{C}(X_{0},X_{1})\otimes \mathcal{C}(X_{1},X_{2})\otimes \cdots \otimes \mathcal{C}(X_{k-1},X_{k})\to \mathcal{C}(X_{0},X_{k}).\] For the sign convention in \(A_\infty\)-categories, we follow \cite{Seidel}, which we note is different from the sign convention used in \cite{NRSSZ}, \cite{CNS}.

Let \(\mathcal{A}\) be a semi-free DGA with generating set \(\mathcal{S}=\mathcal{R}\cup\mathcal{H}\) with a \(P\)-link grading. Let \(\epsilon:\mathcal{A}\to\End(V)\) be a representation of \(\mathcal{A}\) which respects the link grading.

We write \(\boxtimes:=\otimes_{\End(V)}\). Let \(C\) be the free \(\End(V)\)-bimodule generated by \(\mathcal{R}\). We define \[\mathcal{A}^\epsilon:=\bigoplus_{k\geq 0}C^{\boxtimes k};\] and \(\mathcal{A}^\epsilon_+\subset \mathcal{A}^\epsilon\) by \[\mathcal{A}^\epsilon_+:=\bigoplus_{k\geq 1}C^{\boxtimes k}.\]

Let \(\phi_\epsilon:\mathcal{A}\to \mathcal{A}^\epsilon\) be an \(\F\)-algebra map defined by its action on generators: for \(q\in\mathcal{R}\), \(\phi_\epsilon(q):=q+\epsilon(q)\), and for \(\gamma\in R\), \(\phi_\epsilon(\gamma):=\epsilon(\gamma)\). Define \(\partial_\epsilon: \mathcal{A}^\epsilon\to \mathcal{A}^\epsilon\) as follows: for \(q\in\mathcal{R}\), \(\partial_\epsilon(q):=\phi_\epsilon(\bar{\partial}q)\), and for \(A\in\End(V)\), \(\partial_\epsilon(A):=0\). We have that \(\partial_\epsilon\) is an \textit{augmented differential}, that is, \(\partial_\epsilon(\mathcal{A}^\epsilon)\subset \mathcal{A}^\epsilon_+\).

Let \(C^\vee:=\Hom_{\End(V)-\End(V)}(C,\End(V))\) be the dual of \(C\). \(C^\vee\) is the free \(\End(V)\)-bimodule generated by \(\{q^\vee:q\in\mathcal{R}\}\), where we define the grading by \(|q^\vee|=|q|+1\). We define \(A_\infty\)-multiplication \[m_k(\epsilon):(C^\vee)^{\otimes k}\to C^\vee\]
by
\[
m_k(\epsilon)(B_k q_{i_k}^\vee,...,B_1 q_{i_1}^\vee)=(-1)^\sigma\sum_{q\in\mathcal{R},A_i\in\End(V)} \mathrm{Coeff}_{A_1q_{i_1}...A_kq_{i_k}A_{k+1}}(\partial_\epsilon q)\cdot A_1B_1\ldots A_k B_k A_{k+1}q^\vee,
\]
where \(\mathrm{Coeff}\) means coefficient, and \[\sigma=|q_{i_1}^\vee|+|q_{i_2}^\vee|+\cdots-k.\] We will simply write \(m_k\) for \(m_k(\epsilon)\), when there is no confusion.

\begin{prop}
	The \(A_\infty\)-relations are satisfied by the above operations.
\end{prop}
\begin{proof}
	The proof follows exactly as in \cite{CDGG}*{Theorem 4.11}. One caveat is that our signs differ from those in \cite{CDGG}, since the order of composition in our category is the opposite of theirs, which introduces the sign \(\sigma\), see \cite{Seidel}*{Section (1a)}.
\end{proof}

We have the decomposition \(C^\vee=\bigoplus_{i,j\in P}C_{ij}^\vee\), where \(C_{ij}^\vee\) is the free \(\End(V)\)-bimodule generated by \(\{q_{ij}^\vee:q_{ij}\in \mathcal{R}_{ij}\}\).

Assume that \(\epsilon\) is diagonal (this is the case we are particularly interested in, since it corresponds to the case where we have \(k\) generic disjoint Legendrians, each equipped with a representation). Then, we have that the only non-trivial components of \(m_k\) are the composable ones:
\[
m_k:C_{i_0i_1}^\vee\otimes\cdots\otimes C_{i_{k-1} i_k}^\vee \to C_{i_0 i_k}^\vee.
\]

Let \(e_i:=\id\in\End(V_i)\). Define the subspace \(D_{ij}^\vee\subset C_{ij}^\vee\) by \(D_{ij}^\vee:=e_i C_{ij}^\vee e_j\). By construction, \(m_k\) restricts to an \(A_\infty\)-operation
\begin{equation}\label{D A-infinity operation}
m_k:D_{i_0i_1}^\vee\otimes\cdots\otimes D_{i_{k-1} i_k}^\vee \to D_{i_0 i_k}^\vee,
\end{equation}
giving rise to an \(A_\infty\)-category whose object set is \(P\).

The following lemma is useful for showing that the \(A_\infty\)-operation in the partially wrapped contact Fukaya category is well-defined. Let \(P'\subset P\) be a subset. The DGA \(\mathcal{A}_{P'}\) (see Definition \ref{definition quotient DGA}) has an inherited \(P'\)-link grading. We also have that the diagonal representation \(\epsilon\) induces a diagonal representation \(\epsilon':\mathcal{A}_{P'}\to\End(\bigoplus_{i\in P'}V_i)\). The same construction as above gives an \(A_\infty\)-category structure on \(P'\).

\begin{lem}\label{subset of P}
	The \(A_\infty\)-category structure on \(P'\) given by \(\epsilon'\) and the \(A_\infty\)-structure on \(P'\) as a subcategory of \(P\) are identical.
\end{lem}

\subsection{Functoriality}\label{algebraic functoriality}
We discuss the functoriality of the above construction, following the end of \cite{CNS}*{Section 2.1}, which will be used to prove the isotopy invariance of bilinearized contact cohomology, and the isotopy invariance of the partially wrapped contact Fukaya category.

Let \(\mathcal{A},\mathcal{B}\) be semi-free DGAs over the same coefficient ring \(R\) (\(R\) being generated by the same generating set \(\mathcal{H}\)), equipped with \(P\)-link gradings that take equal values on \(\mathcal{H}\). Let \(f:\mathcal{A}\to\mathcal{B}\) be a DGA morphism that respects link grading, and is identity when restricted to \(\mathcal{H}\). Let \(\epsilon':\mathcal{B}\to\End(V)\) be a representation respecting link grading. By composing with \(f\), we obtain a representation \(\epsilon:=\epsilon'\circ f:\mathcal{A}\to\End(V)\) that respects the link grading. To distinguish between the \(C\)'s, we write \[\mathcal{A}^\epsilon=\bigoplus_{k\geq0}C_\mathcal{A}^{\boxtimes k}, \mathcal{B}^\epsilon=\bigoplus_{k\geq0}C_\mathcal{B}^{\boxtimes k}.\]

We define an \(\F\)-algebra morphism \(f^\epsilon:\mathcal{A}^{\epsilon'}\to \mathcal{B}^\epsilon\) by its action on generators:
\begin{align*}
	&f^\epsilon(q)=\phi_{\epsilon'}f(q)-\epsilon(q) && q\in\mathcal{R}&\\
	&f^\epsilon(\alpha)=\alpha && \alpha\in\End(V).&
\end{align*}
One can check on the generators that \(f^\epsilon\circ \partial_{\epsilon'}=\partial_\epsilon\circ f^\epsilon\). We have that the constant term in \(f^\epsilon(q)\) is \(\epsilon'f(q)-\epsilon(q)=0\), so \(f^\epsilon (C)\subset \mathcal{B}_+^{\epsilon'}\). So we can write \(f^\epsilon|_{C}=f_1\oplus f_2\oplus\cdots\) where each \(f_k\) has image in \(C_\mathcal{B}^{\boxtimes k}\). Dualizing all the \(f_i\) as in the construction of the \(m_k\), we get a family of maps \(\{F_k\}_{k\in\N}\) which forms an \(A_\infty\)-algebra morphism \(F:(C_\mathcal{B}^\vee,\{m_k(\epsilon')\})\to (C_\mathcal{A}^\vee,\{m_k(\epsilon)\})\). When we assume that the representation \(\epsilon'\) is diagonal, \(F\) can be upgraded to an \(A_\infty\)-functor between the categories defined in \ref{D A-infinity operation}. If \(f\) is a stable tame isomorphism, then on the level of cohomology, the functor \(F\) is an isomorphism on the morphisms.

The following lemma is useful for showing that when proving the invariance of the partially wrapped contact Fukaya category, the functor we construct is well-defined.

\begin{lem}\label{functor subset of P}
	Let \(P'\subset P\) be a subset. Lemma \ref{subset of P} gives two \(A_\infty\)-category structures on \(P'\), one from \(\mathcal{A}'\) and one from \(\mathcal{B}'\). The functor \(F'\) between them constructed in this section using the induced map \(f':\mathcal{A}'\to\mathcal{B}'\) is identical to the restricted functor \(F_{|P'}\).
\end{lem}

\section{Partially wrapped contact Fukaya categories}\label{section partially wrapped contact Fukaya categories}

In this section, we construct the partially wrapped contact Fukaya categories of \(1\)-jet spaces, and prove various structural results. We first defined the morphisms and \(A_\infty\)-operations in the unwrapped setting. Then, we use localization at continuation elements to define the partially wrapped contact Fukaya categories, and prove their isotopy invariance. Along the way, we give a characterization of the representation category of a Legendrian as a full subcategory of a partially wrapped contact Fukaya category. To facilitate the computation of partially wrapped contact Fukaya categories in Section \ref{section equivalence}, we also characterize the partially wrapped contact Fukaya categories with empty stops, and prove a wrapping exact triangle.

\subsection{The unwrapped contact Fukaya pre-category}\label{section unwrapped}

We construct the unwrapped contact Fukaya pre-category \(\mathcal{F}uk(V,J)\) of a 1-jet space \(V=J^1M\) with a generic compatible almost complex structure \(J\) on \(V\times\R\), which is an \(A_\infty\)-pre-category in the sense of \cite{GPS2}*{Definition 2.6}, where almost all morphisms and compositions are defined. More work is needed to define all morphisms and compositions to make this an actual \(A_\infty\)-category, but we do not need that in this paper.

\begin{rem}
   The construction of the unwrapped contact Fukaya pre-category and the partially wrapped contact Fukaya category in the more general case \(V=P\times \R\) follows as in this paper with virtually no change, except we need to pick paths to specify the gradings of Reeb chords, instead of simply picking a Maslov potential. Rather than working with Legendrians with vanishing Maslov class, we can pick a non-negative integer \(m\) and work with Legendrians with Maslov class divisble by \(m\), in which case all chain complexes would be \(\Z/2m\Z\)-graded. We can also drop the relatively pin assumption, if we work with characteristic 2 coefficients.
\end{rem}

Fix a 1-jet space \(V=J^1M\), with Reeb vector field \(\frac{\partial}{\partial z}\), a field \(\F\). We say that a compact relatively pin Legendrian \(\Lambda\subset V\) is compatible with the almost complex structure \(J\), if the Chekanov-Eliashberg DGA of \(\Lambda\) can be defined using \(J\). Objects of \(\cfuk(J^1M)\) are Legendrian \(R\)-branes:

\begin{defn}\label{definition R-brane}
A \emph{Legendrian \(R\)-brane} (\(R\) is for representation) consists of the data of a compact Legendrian submanifold \(\Lambda\subset J^1M\) compatible with \(J\), relatively pin with repect to the class \(w_2(M)\), with vanishing Maslov class, equipped with a Maslov potential, and a representation \(\epsilon\) of \(\Lambda\) over \(\F\). Additionally, every component \(\Lambda_i\) of \(\Lambda\) is equipped with a basepoint \(p_i\). Every point on \(\Lambda_i\) is equipped with a path on \(\Lambda_i\) to \(p_i\) (the endpoint path). By abuse of notation, we will often denote a Legendrian \(R\)-brane by its underlying Legendrian \(\Lambda\), and we will denote its representation \(\epsilon_{\Lambda}:\mathcal{A}(\Lambda)\to \End_{\F}(V^\Lambda)\).

An \emph{isotopy of Legendrian \(R\)-branes} \(\Lambda^t:t\in[0,1]\) is an isotopy of the underlying Legendrians, such that the relative pin structure, Maslov potential, basepoints, endpoint paths are simply carried forward smoothly. The representation \(\epsilon_{\Lambda^1}\) depends on \(\epsilon_{\Lambda^0}\) in the following way: the proof of invariance of the Chekanov-Eliashberg DGA from \cite{EES3} gives a stable tame isomorphism \(\mathcal{A}(\Lambda^1)\to\mathcal{A}(\Lambda^0)\). Composing by this stable tame isomorphism, we get a map from the set of representations of \(\Lambda^0\) to that of \(\Lambda^1\), which sends \(\epsilon_{\Lambda^0}\) to \(\epsilon_{\Lambda^1}\).
\end{defn}

\begin{rem}
   We fixed the background cohomology class \(w_2(M)\in H^2(J^1M;\Z/2)\), but the construction of the partially wrapped contact Fukaya category carries through if we replace this with an arbitrary class \(b\in H^2(J^1M;\Z/2)\). The wrapping exact triangle would still hold. The difference is that we make use of sections of \(J^1M\), which are guaranteed to be relatively pin with respect to \(w_2(M)\), but might not be relatively pin with respect to some other class.
\end{rem}

Let \(\Lambda_0,\Lambda_1\) be Legendrian \(R\)-branes. Assume that \(\Lambda_0\) and \(\Lambda_1\) are disjoint, and \(\Lambda_0\cup\Lambda_1\) is chord generic. Let \(R_{\Lambda_0,\Lambda_1}\) denote the set of Reeb chords from \(\Lambda_0\) to \(\Lambda_1\). Then, we define the morphism space \(Hom_{\cfuk(V,J)}(\Lambda_0,\Lambda_1)\) (denoted \(LCC(\Lambda_0,\Lambda_1)\), called the \emph{bilinearized contact cochain complex}), to be the space \(D_{01}^\vee\) as in Section \ref{section augmented DGAs}. Spelling it out, 
\[
LCC(\Lambda_0,\Lambda_1):=\bigoplus_{b\in R_{\Lambda_0,\Lambda_1}}\Hom_{\F}(V^{\Lambda_0}_{c(b)},V^{\Lambda_1}_{r(b)})b^\vee
\]
is a \(\Z/2m\)-graded vector space, and the grading of a generator \(B b^\vee\in LCC(\Lambda_0,\Lambda_1)\) is 1 less than the grading of \(b\) as a generator of the Chekanov-Eliashberg DGA of \(\Lambda_0\cup \Lambda_1\). The link grading \(c,r\) is given by decomposing \(\Lambda_0\cup\Lambda_1\) into connected components.

\begin{figure}
   \centering
	\includegraphics{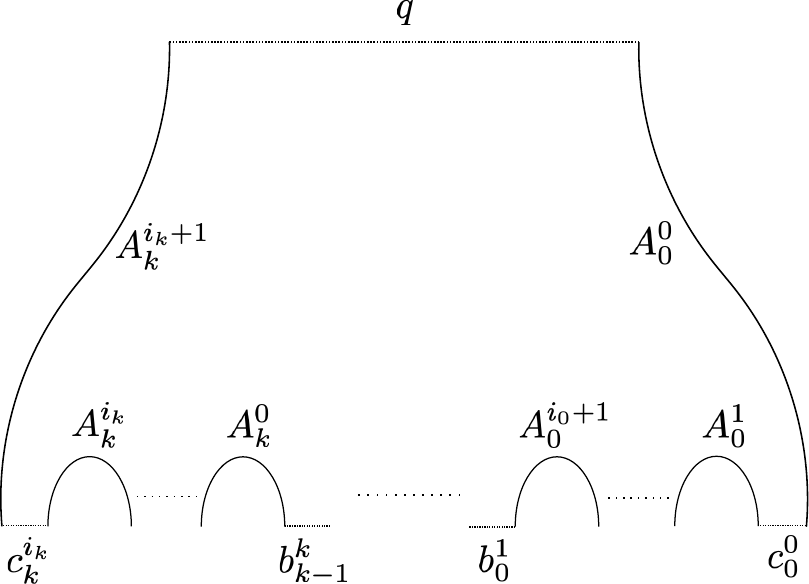}
	\caption{Holomorphic disks counted in the \(A_\infty\)-operations}
	\label{figure A_infinity disk}
\end{figure}

Let \(\Lambda_0,\ldots,\Lambda_k\) be Legendrian \(R\)-branes. Assume that \(\Lambda_0,\ldots,\Lambda_k\) are pairwise disjoint, and \(\bigcup_{0\leq i\leq k}\Lambda_i\) is generic admissible (so that we can count Reeb chords and holomorphic disks on it). For each \(0\leq i\leq k-1\), let \(b_{i,i+1}\) be a Reeb chord from \(\Lambda_i\) to \(\Lambda_{i+1}\), and let \(B_i\in\Hom_{\F}(V^{\Lambda_i}_{c(b_{i,i+1})},V^{\Lambda_{i+1}}_{r(b_{i,i+1})})\). Then, we define the \(A_\infty\)-operation \(m_k\) to be the \(A_\infty\)-operation as in \ref{D A-infinity operation}.

Given a finite subset \(P\in \{0,\ldots,k\}\), we have that \(\mathcal{A}(\bigcup_{i\in P}\Lambda_i)\) is canonically isomorphic to \(\mathcal{A}_P\). Lemma \ref{subset of P} then tells us that the \(A_\infty\)-operation is well-defined, that is, the operations \[m_{|P|}:LCC(\Lambda_{i_0},\Lambda_{i_1})\otimes\cdots\otimes LCC(\Lambda_{i_{|P|-1}},\Lambda_{i_{|P|}})\to LCC(\Lambda_{0},\Lambda_{i_{|P|}})\] for \(i_j\in P\) does not depend on whether it comes from considering \(\bigcup_{0\leq i\leq k}\Lambda_i\) or \(\bigcup_{P}\Lambda_i\).

More concretely, \(m_k\) is given by the following formula, which gives another way to show the well-definedness of \(m_k\):
\begin{align*}
m_k(&B_{0}b_{0,1}^\vee,\ldots,B_{k-1}b_{k-1,k}^\vee)\\
=&\sum (-1)^{\bigstar}\left |\mathcal{M}_{\bar{A}}(q;c^{i_k}_{k},\ldots,c^{0}_{k},b_{k-1,k},c^{i_{k-1}}_{k-1},\ldots,c^{0}_{k-1},b_{k-2,k-1},c^{i_{k-2}}_{k-2}\ldots\ldots,b_{0,1}c^{i_0}_{0},\ldots,c^{0}_{0})\right |\\
&\epsilon_{\Lambda_0}(A^{0}_{0}c^{0}_{0}\cdots A^{i_{0}}_{0}c^{i_{0}}_{0}A_{0}^{i_0+1})B_{0}\epsilon_{\Lambda_1}(A^{0}_{1}c^{0}_{1}\cdots A^{i_{1}}_{1}c^{i_{1}}_{1}A_{1}^{i_1+1})B_{1}\cdots B_{k-1}\epsilon_{\Lambda_k}(A^{0}_{k}c^{0}_{k}\cdots A^{i_{k}}_{k}c^{i_{k}}_{k}A_{k}^{i_k+1})q^\vee,
\end{align*}
where the sum is over all tuples of non-negative integers \((i_0,\ldots,i_k)\), over Reeb chords \(c^{i}_{j}\) from \(\Lambda_j\) to \(\Lambda_j\), and all possible \(A^{i}_{j}\in \pi_1(\Lambda_j)\). The sign is \(\bigstar=|b_{01}^\vee|+|b_{12}^\vee|+\cdots-k+(n-1)(|q^\vee|)\). See Figure \ref{figure A_infinity disk} for the disks that correspond to each summand.

\subsection{Bilinearized contact cohomology and the continuation element}\label{section continuation elements}

We expand on bilinearized contact cohomology, and define the continuation element. See \cite{GPS2}*{Section 2.1} for an account of stops and wrappings of Legendrians. We will freely use the terminologies of \cite{GPS2}*{Section 2.1} in the rest of this paper. The story here is roughly analogues to \cite{GPS1}*{Section 3.3}.

Let \(\Lambda_0,\Lambda_2\) be Legendrian \(R\)-branes. The \emph{bilinearized contact cohomology}, \(LCH(\Lambda_0,\Lambda_1)\), is the cohomology of the bilinearized contact cochain complex \((LCC(\Lambda_0,\Lambda_1),m_1)\). Note that this is different from what is called bilinearized contact cohomology in \cite{BC}.

\begin{prop}\label{invariance of BLCC}
	Bilinearized Legendrian contact cohomology is invariant under Legendrian isotopies. That is, let \(H\) be the time-\(1\) map of a contact isotopy, then \[LCH(\Lambda,\Lambda')\cong LCH(H(\Lambda),H(\Lambda')).\]
\end{prop}
\begin{proof}
	This is an application of the stable tame isomorphism invariance of the Chekanov-Eliashberg DGA, Theorem \ref{EES theorem}.
\end{proof}

We write \[LCH(\Lambda):=LCH(\Lambda,\Lambda^+),\] where \(\Lambda^+\) is a generic small positive pushoff of \(\Lambda\) (small enough so that on \(\Lambda\cup\Lambda^+\), rigid holomorphic disks on correspond to rigid quantum flow trees, more precisely, small enough so that \cite{Liu25}*{Proposition 3.13} applies). Reeb chords from \(\Lambda\) to \(\Lambda^+\) consists of long Reeb chords, which are in bijection with self Reeb chords of \(\Lambda\), and short Reeb chords, which are in bijection with the set of critical points of the Morse function \(f:\Lambda\to\R\) used to pushoff \(\Lambda\) to obtain \(\Lambda^+\). Write \(\min_i\) for the short Reeb chords corresponding to index 0 critical points of \(f\). Write 
\[
	c^\vee:=(-1)^{\dagger}\left[\sum I{\min}_i^\vee\right],
\]
where \(I\) is the identity map on \(V^\Lambda\), and \(\dagger\) is defined after Lemma \ref{min min curve}.

\begin{prop}\label{unitality of c}
	\(LCH(\Lambda)\) is a unital algebra, with unit given by \(c^\vee\). Furthermore, the left module \(LCH(\Lambda,\Lambda')\) and the right module \(LCH(\Lambda',\Lambda)\) over \(LCH(\Lambda)\) are unital.
\end{prop}
\begin{proof}
	See Appendix \ref{sign appendix}.
\end{proof}

\begin{defn}
	Given any positive isotopy \(\Lambda^t\), we define a \emph{continuation element} \[c^\vee(\Lambda_t)\in LCH(\Lambda^0,\Lambda^1)\] as follows. For sufficiently small positive isotopies (small enough so that \cite{Liu25}*{Proposition 3.13} applies to the \(2\)-copy \(\Lambda^0\cup\Lambda^1\)), \(c^\vee(\Lambda^t)\) is simply defined as the unit \(c^\vee\) in \(LCH(\Lambda)\). For positive isotopies that is contained in a Weinstein neighborhood of \(\Lambda^0\), we isotope the 2-copy Legendrian \(\Lambda^0\cup\Lambda^1\) to be close enough so that \cite{Liu25}*{Proposition 3.13} applies, and then use Proposition \ref{invariance of BLCC} to obtain \(c^\vee(\Lambda^t)\). For an arbitrary isotopy, break it up into small isotopies such that each small isotopy happens within a Weinstein neighborhood of its starting point, and define \(c^\vee(\Lambda^t)\) as the composition of the continuation elements of each of the pieces.
\end{defn}

Note that for the case of an isotopy contained in a Weinstein neighborhood, \(c^\vee(\Lambda^t)\) is well-defined, since given any 2 isotopies from \(\Lambda^0\cup(\Lambda^0)^+\) (where \((\Lambda^0)^+\) is a positive pushoff that is small enough for \cite{Liu25}*{Proposition 3.13}) to \(\Lambda^0\cup\Lambda^1\), we can concatenate one with the inverse of the other, and obtain a loop of Legendrians based at \(\Lambda^0\cup(\Lambda^0)^+\). By Lemma \ref{isotopy sends c to c}, any such loop induces a map that sends \(c^\vee\in LCH(\Lambda)\) to itself. Therefore, the 2 isotopies induce maps on \(LCH\) that send \(c^\vee\in LCH(\Lambda)\) to the same element in \(LCH(\Lambda^0,\Lambda^1)\).

\begin{lem}\label{continuation element lemma}
	The continuation element satisfies the following properties:
	\begin{enumerate}
		\item \(c^\vee(\Lambda^t)\) is independent of how the isotopy is broken up into small isotopies.
		\item \(c^\vee(\Lambda^t)\) is invariant under deformation of the isotopy \(\Lambda_t\) fixed at endpoints.
		\item A composition of continuation elements coincides with the continuation
element of the concatenation of isotopies.
	\end{enumerate}
\end{lem}

\begin{proof}
	For \((1)\): let \(\Lambda^t,t\in[0,1]\) be a positive isotopy. Given \(0\leq s<t\leq 1\) for \(t-s\) small enough, let \(c_{s,t}^\vee\in LCH(\Lambda^s,\Lambda^t)\) denote the continuation element. Since \(c_{s,t}^\vee\) is the unit element in \(LCH(\Lambda^s)\), we have that \(c_{x,y}^\vee\cdot c_{y,z}^\vee=c_{x,z}^\vee\), if \(x<y<z\) are close enough. This proves that \(c^\vee(\Lambda^t)\) is independent of the choice of a sufficiently fine partition \(0=t_0<\cdots <t_n=1\). This also proves \((3)\).

	For \((2)\): we can break up any deformation fixed at endpoints into small deformations fixed at endpoints. Let \(\Lambda^{t,s}\) be such a small deformation fixed at endpoints. We can find partitions \(0=t_0<\cdots <t_n=1\) and \(0=t'_0<\cdots <t'_n=1\) such that for each \(0< i < n\), there is a small positive isotopy from \(\Lambda^{{t_i},0}\) to \(\Lambda^{{t'_i},0}\). The proof then follows from the equation \(c_{x,y}^\vee\cdot c_{y,z}^\vee=c_{x,z}^\vee\).
\end{proof}

\subsection{Partially wrapped contact Fukaya categories}

Let \(V=J^1M\) be a 1-jet space, and let \(\mathfrak{f}\subset V\) be a closed subset (called a \emph{stop}). We construct the partially wrapped contact Fukaya category \(\mathcal{V}(V,\mathfrak{f})\), which is a cohomologically unital \(A_\infty\)-category, and prove that it is well-defined. This section is analogous to the construction of wrapped Fukaya categories in \cite{GPS1}*{Section 3.5}.

Let \(J\) be a generic compatible almost complex structure on \(V\times\R\). Fix a countable collection \(\tilde{I}\) of Legendrian submanifolds equipped with all the data of \(R\)-branes except for representations, such that for each isotopy class of compact Legendrian submanifolds of \(V\setminus\mathfrak{f}\) relatively pin with respect to \(w_2(M)\), with vanishing Maslov class, there is a representative \(\tilde{\Lambda}\in \tilde{I}\).

\begin{defn}
	(Poset \(\mathcal{O}\)) For each Legendrian \(\tilde{\Lambda}\in \tilde{I}\), choose a cofinal wrapping sequence
	\[
	\tilde{\Lambda}=\tilde{\Lambda}^{(0)}\rightsquigarrow \tilde{\Lambda}^{(1)}\rightsquigarrow \cdots
	\]
	of morphisms in the wrapping category \((\Lambda\rightsquigarrow -)^+_{V\setminus\mathfrak{f}}\). Define \(\mathcal{O}\) as \[\bigcup_{\tilde{\Lambda}^{(i)}\in \Z_{\geq0}\times\tilde{I}}\{\Lambda^{(i)}: \Lambda^{(i)} \text{ is a Legendrian \(R\)-brane given by equipping } \tilde{\Lambda}^{(i)} \text{with a representation}\}.\]

	We equip \(\mathcal{O}\) with the partial order inherited from the order on \(\Z_{\geq0}\), namely \(\Lambda^{(i)}<\Lambda'^{(i')}\) iff \(i<i'\). By choosing the \(\tilde{\Lambda}^{(i)}\) generically, we ensure that the union any finite totally ordered collection of Legendrians appearing in \(\mathcal{O}\) is generic admissible.
\end{defn}

\begin{defn}\label{definition of category O}
	(Category \(\mathcal{O}\)) We turn \(\mathcal{O}\) into a strictly unital \(A_\infty\)-category with the following morphism spaces:
	\[
	\mathcal{O}(\Lambda,\Lambda'):=
	\begin{cases}
		LCC(\Lambda,\Lambda') & \Lambda < \Lambda'\\
		\F[0] & \Lambda=\Lambda'\\
		0 &\text{otherwise.}
	\end{cases}
	\]
	The \(A_\infty\)-operation
	\[
	m_k:\mathcal{O}(\Lambda_0,\Lambda_1)\otimes\cdots\otimes \mathcal{O}(\Lambda_{k-1},\Lambda_k)\to \mathcal{O}(\Lambda_0,\Lambda_k)
	\]
	is defined to be the same as the \(A_\infty\)-operation in the unwrapped contact Fukaya pre-category \(\cfuk(V,J)\), if \(\Lambda_0<\cdots<\Lambda_k\). The only other possibly nontrivial operations, when some \(\Lambda_i=\Lambda_{i+1}\), are formally fixed by strict unitality.
\end{defn}

\begin{rem}
	Note that although the partial wrapping can be viewed as the flow of a cutoff Reeb field given by the stop \(\mathfrak{f}\), \(LCC(\Lambda,\Lambda')\) is generated by chords of the vertical Reeb field \(\frac{\partial}{\partial z}\).
\end{rem}

\begin{defn}\label{definition of partially wrapped contact Fukaya categories}
	Denote by \(C\) the set of all continuation elements \(c\in LCH(\Lambda^{(i)},\Lambda^{(i+1)})\), for all \(\Lambda^{(i)}\in\mathcal{O}\). We define the \emph{partially wrapped contact Fukaya category} \(\mathcal{V}(V,\mathfrak{f})\) to be the localized category \(\mathcal{V}(V,\mathfrak{f}):=\mathcal{O}[C^{-1}]\).
\end{defn}


\begin{lem}\label{lemma 2 isomorphisms}
	The natural maps 
	\[
	\varinjlim_i H^*\mathcal{O}(\Lambda,\Lambda'^{(i)})\xrightarrow\sim\varinjlim_i H^*\mathcal{V}(\Lambda,\Lambda'^{(i)})\xleftarrow\sim H^*\mathcal{V}(\Lambda,\Lambda')
	\]
	are both isomorphisms.
\end{lem}
\begin{proof}
	The proof is similar to the proof of \cite{GPS1}*{Lemma 3.37}. Since the continuation maps are isomorphisms in \(\mathcal{V}\), all maps in the colimit in the middle term are isomorphisms. Therefore, the rightmost arrow is an isomorphism.

   By \cite{GPS1}*{Lemma 3.16}, to show that the leftmost arrow is an isomorphism, it is enough to show that \(\varinjlim_i \Lambda'^{(i)}\) is right \(C\)-local, meaning every continuation map \(\Lambda''^{(j)}\to\Lambda''^{(j+1)}\) induces an isomorphism \[\varinjlim_i LCH(\Lambda''^{(j)},\Lambda'^{(i)})\to \varinjlim_i LCH(\Lambda''^{(j+1)},\Lambda'^{(i)})\] is an isomorphism. This map is simply \(H^*\mathcal{V}(\Lambda''^{(j)},\Lambda')\to H^*\mathcal{V}(\Lambda''^{(j+1)},\Lambda')\), which is an isomorphism since it is given by composition with an isomorphism.
\end{proof}


To prove that \(\mathcal{V}(V,\mathfrak{f})\) is well-defined, we give a more flexible construction of \(\mathcal{V}\), as follows. Let \(\mathcal{O}\) be any poset of Legendrian \(R\)-branes, such that: the set of underlying Legendrians of elements in \(\mathcal{O}\) is countable, and contain every Legendrian isotopy class; for every finite totally ordered collection of Legendrian \(R\)-branes in \(\mathcal{O}\), their union is generic admissible. Let \(C\) be any collection of elements of \(LCH(\Lambda,\Lambda')\) for \(\Lambda < \Lambda'\in\mathcal{O}\) consisting only of continuation elements for various positive isotopies from \(\Lambda\) to \(\Lambda'\). Now assume that the pair \((\mathcal{O},C)\) satisfies the following cofinality property:
\begin{itemize}
	\item For every \(\Lambda\in\mathcal{O}\), there is a cofinal sequence \(\Lambda=\Lambda^{(0)}<\Lambda^{(1)}<\cdots\) in \(\mathcal{O}\) together with positive isotopies \(\Lambda^{(i)}\rightsquigarrow\Lambda^{(i+1)}\) such that \(\Lambda=\Lambda^{(0)}\rightsquigarrow\Lambda^{(1)}\rightsquigarrow\cdots\) is cofinal in the wrapping category \((\Lambda\rightsquigarrow -)^+_{V\setminus\mathfrak{f}}\), and each of the resulting continuation elements in \(LCH(\Lambda^{(i)},\Lambda^{(i+1)})\) is in \(C\).
\end{itemize}

We can then equip it with the structure of an \(A_\infty\)-category like in Definition \ref{definition of category O}, and localize to get \(\mathcal{V}(V,\mathfrak{f}):=\mathcal{O}[C^{-1}]\) like in Definition \ref{definition of partially wrapped contact Fukaya categories}.

\begin{prop}\label{proposition well-defined}
	Up to quasi-equivalence, the category \(\mathcal{V}(V,\mathfrak{f})\) defined above is independent of the choices of \(\mathcal{O},C\).
\end{prop}
Note that in this proposition, we do not address the issue of well-definedness with respect to different almost complex structures. We will address that in Theorem \ref{isotopy invariance}.
\begin{proof}
	The proof goes exactly as \cite{GPS1}*{Proposition 3.39}. If \(\mathcal{O}\subseteq\mathcal{O}'\) and \(C\subseteq C'\), then there is a natural functor \(\mathcal{V}\to\mathcal{V}'\). By Lemma \ref{lemma 2 isomorphisms}, this functor is a quasi-equivalence.

   It thus suffices to show that any \((\mathcal{O},C),(\mathcal{O}',C')\) can be included into a third. Consider \(\mathcal{O}\sqcup\mathcal{O}'\) with no order relations between \(\mathcal{O}\) and \(\mathcal{O}'\). Let \(\mathcal{O}'':=\Z_{\geq0}\times(\mathcal{O}\sqcup\mathcal{O}')\) with its lexicographical partial order, choosing for every Legendrian \(R\)-brane in \(\mathcal{O}\sqcup\mathcal{O}'\) a generic cofinal wrapping sequence such that: for any two Legendrian \(R\)-branes \(\Lambda,\Lambda'\in \mathcal{O}\sqcup\mathcal{O}'\) that have the same underlying Legendrian, the underlying Legendrians of their cofinal wrapping sequences are chosen to be the same. Let \(C''\) denote the union of the resulting continuation elements together with \(C,C'\). Then, \((\mathcal{O}'',C'')\) satisfies the cofinality property, and \((\mathcal{O},C),(\mathcal{O}',C')\subseteq (\mathcal{O}'',C'')\).
\end{proof}

\subsection{The representation category of a Legendrian}

The augmentation category \(\mathcal{A}ug_+(\Lambda)\) of a Legendrian \(\Lambda\in V\) (defined for Legendrians in \(\R^3\) \cite{NRSSZ}, and more generally \cite{Liu25}), and its higher rank generalization, the rank \(n\) representation category \(\mathcal{R}ep_n(\Lambda)\) (defined for Legendrians in \(\R^3\) \cite{CNS}, and can be defined more generally using the localization methods of \cite{Liu25}) are naturally full subcategories of the finite rank representation category \(\mathcal{R}ep(\Lambda)\), which is naturally a full subcategory of the partially wrapped contact Fukaya category, as follows:

\begin{defn}
	Let \(\Lambda\in V\) be a Legendrian equipped with all data of an \(R\)-brane except for a representation. Let \(\Lambda^+\) be a small positive pushoff of \(\Lambda\). Then, we define the \textit{representation category} \(\mathcal{R}ep_+(\Lambda)\) to be the full subcategory of \(\mathcal{V}(V,\Lambda^+)\) consisting of Legendrian \(R\)-branes whose underlying Legendrian is \(\Lambda\). Write \(\mathcal{R}ep_+^{n}(\Lambda)\) for the subcategory of \(\mathcal{R}ep_+(\Lambda)\) where the representation on the Legendrian \(R\)-brane is \(n\)-dimensional. Write \(\mathcal{A}ug_+(\Lambda)\) for \(\mathcal{R}ep_1(\Lambda)\).
\end{defn}

\begin{rem}
	To see that \(\mathcal{R}ep_+(\Lambda)\) as defined above is the same as the rank \(n\) representation category defined using small positive pushoffs in the literature, note that in the definition of \(\mathcal{R}ep_+(\Lambda)\) as above, we can take a cofinal wrapping sequence of \(\Lambda\) that simply wraps into \(\Lambda^+\), which consist of exactly the small positive pushoffs as in the definition of the augmentation category and the representation category in the literature.
\end{rem}

\begin{rem}
   The above construction of \(\mathcal{R}ep\) generalizes previous work \cites{NRSSZ,CNS,Liu25} in several ways. It generalizes the construction of the representation category \(\mathcal{R}ep_n\) in \(\R^3\) in \cite{CNS} (which already generalizes the construction of the augmentation category in \(\R^3\) in \cite{NRSSZ} to arbitrary ranks) to arbitrary \(1\)-jet spaces \(J^1M\) (or even more generally, \(P\times\R\)), and with all finite rank representations in the same category rather than a category for each rank. It generalizes the construction of the augmentation category \(\mathcal{A}ug_+\) in arbitrary \(P\times\R\) in \cite{Liu25} from rank 1 to arbitrary rank. Since we deal with the issue of orientations of moduli spaces in this paper, this also generalizes the characteristic \(2\) coefficient fields in \cite{Liu25}. Note that one can also construct \(\mathcal{R}ep\) using methods in \cite{Liu25} directly, together with the discussion of orientations of moduli spaces in this paper.
\end{rem}

The reason we define \(\mathcal{R}ep_+(\Lambda)\) here, right after defining partially wrapped contact Fukaya categories, is that the invariance of \(\mathcal{R}ep_+(\Lambda)\) will be proved and used in the proof of the invariance of partially wrapped contact Fukaya categories, Theorem \ref{isotopy invariance}.

\begin{thm}
	\(\mathcal{R}ep_+(\Lambda)\) is invariant under Legendrian isotopy of \(\Lambda\).
\end{thm}
\begin{proof}
	This theorem will be proved as a step in the proof of Theorem \ref{isotopy invariance}.
\end{proof}

\begin{rem}
   In \cite{NRSSZ}, it was proved that for any Legendrian link \(\Lambda\subset\R^3\), there is an equivalence between the augmentation category and the category of sheaves on \(\R^2\) with microlocal rank \(1\) with microsupport on \(\Lambda\). It was conjectured in \cite{CNS} that for any Legendrian \(\Lambda\subset J^1M\) in any \(1\)-jet space, a similar equivalence holds between \(\mathcal{R}ep(\Lambda)\) and the category of sheaves on \(M\times\R\) with microsupport in \(\Lambda\). This conjecture is false, by \cites{Gao,GL}.
\end{rem}

\subsection{Isotopy invariance}\label{section invariance}

\begin{thm}\label{isotopy invariance}
	The partially wrapped contact Fukaya category is invariant under a contact isotopy of the stop, and is independent of the almost complex structure. That is, if \(\mathfrak{f},\mathfrak{f}'\) are contact isotopic, then \(\mathcal{V}(V,\mathfrak{f})\) and \(\mathcal{V}(V,\mathfrak{f}')\) (defined using possibly different almost complex structures), are quasi-equivalent.
\end{thm}
\begin{proof}
	We first define a functor \(F':\mathcal{O}'\to\mathcal{O}\). For some choice of \(\mathcal{O}'\) that depends on the choice of \(\mathcal{O}\). We will then prove that \(F'\) maps the continuation elements to the continuation elements, thus inducing a functor \(F\) from \(\mathcal{V}(V,\mathfrak{f})\) to \(\mathcal{V}(V,\mathfrak{f}')\), which we will show to be a quasi-equivalence.

	Let \(H\) be the time-1 map of a contact isotopy of \(V\), with \(H(\mathfrak{f})=\mathfrak{f}'\). Given the choices involved in defining \(\mathcal{O}\), that is, an almost complex structure \(J\) and the choice of the poset \(\mathcal{O}\), there is a poset \(\mathcal{O}':=H(\mathcal{O})\) whose Legendrians are all disjoint from \(\mathfrak{f}\). As long as one picks \(H\) and \(\mathcal{O}\) generically, there exists some almost complex structure \(J'\) such that for every finite totally ordered collection of Legendrian \(R\)-branes in \(\mathcal{O}'\), their union is generic admissible. We define an \(A_\infty\)-functor from \(\mathcal{O}'\) to \(\mathcal{O}\) as follows: on the level of objects, send each \(H(\Lambda)\) to \(\Lambda\). For each finite totally ordered collection \(P\) of Legendrian \(R\)-branes in \(\mathcal{O}'\), the contact isotopy \(H\) gives rise to a stable tame isomorphism DGA map respecting link grading
	\[
	f:\mathcal{A}\left(\bigcup_{i\in P}\Lambda_{i}\right)\to \mathcal{B}\left(\bigcup_{i\in P}H(\Lambda_{i})\right),
	\]
	which by a restriction of the construction in Section \ref{algebraic functoriality}, gives rise to a unital \(A_\infty\)-functor \(F'_{|H(P)}\) from the subcategory \(H(P)\in \mathcal{O}'\) to \(P\in \mathcal{O}\) (on the endomorphism spaces, \(F'_{|H(P)}\) is formally fixed by unitality). By Lemma \ref{functor subset of P}, we can patch these \(A_\infty\)-functors together to an \(A_\infty\)-functor \(F':\mathcal{O}'\to\mathcal{O}\). Since each \(F'_{|H(P)}\) is induced from a stable tame isomorphism \(f\), we have that on the level of cohomology, \(H^*F'\) is faithful, and sometimes full: \[H^*F:H^*\mathcal{V}(\Lambda^{(i)},\Lambda'^{(j)})\to H^*\mathcal{V}(F(\Lambda^{(i)}),F(\Lambda'^{(j)}))\] is an isomorphism if \(i<j\), and is the zero map if \(i\ge j\).

	By Lemma \ref{isotopy sends c to c}, \(F'\) sends the continuation morphisms to the continuation morphisms. So, \(F'\) descends to a functor between the localized categories \(F:\mathcal{V}(V,\mathfrak{f}')\to\mathcal{V}(V,\mathfrak{f})\).

	To show that \(F\) is a quasi-equivalence, we need to check that \(H^*F\) is fully faithful and essentially surjective. Since \(F'\) is almost fully faithful, and the left most term in Lemma \ref{lemma 2 isomorphisms} is a colimit where the problematic morphism spaces go away (since for high enough \(i\), \(H^*F'\) on \(H^*\mathcal{O}'(\Lambda,\Lambda'^{(i)})\) is an isomorphism), we can use the isomorphism between the left most term and the right most term in Lemma \ref{lemma 2 isomorphisms} to conclude that \(H^*F\) is an isomorphism on \(H^*\mathcal{V}(\Lambda,\Lambda')\). This shows that \(F\) is fully faithful.

	At this point we take a detour to prove the invariance of \(\mathcal{R}ep_+(\Lambda)\). Clearly, \(F\) restricts to a functor \(\mathcal{R}ep_+(H(\Lambda))\to \mathcal{R}ep_+(\Lambda)\), which is fully faithful. To show that this restriction is a quasi-equivalence, we note that the map on objects of \(\mathcal{R}ep_+(\Lambda)\) is given by pre-composition with a stable tame isomorphism \(\mathcal{A}(\Lambda)\to\mathcal{A}(H(\Lambda))\). We break up the isotopy \(H\) into small pieces so that we break up this stable tame isomorphism into a composition of stabilizations, destabilisations, and tame isomorphisms. For tame isomorphisms, pre-composing gives a bijection map on the set of representations, therefore the functor is essentially surjective, and therefore a quasi-equivalence. For stabilization, pre-composing gives a surjective map on the set of representations, therefore the functor is essentially surjective, and therefore a quasi-equivalence. In the case of destabilization, we use the fact that it is the inverse of stabilization, and the fact that quasi-equivalence is an equivalence relation. This proves the invariance of \(\mathcal{R}ep_+(\Lambda)\).

	Now we prove the essential surjectivity of \(F:\mathcal{V}(V,\mathfrak{f}')\to\mathcal{V}(V,\mathfrak{f})\). Let \(\Lambda\subset \mathcal{O}\), and let \(\Lambda^+\) be a small positive pushoff. Consider the categories \(\mathcal{V}(V,\mathfrak{f}\cup \Lambda^+)\) and \(\mathcal{V}(V,\mathfrak{f}'\cup H(\Lambda^+))\), which contain \(\mathcal{R}ep_+(\Lambda)\) and \(\mathcal{R}ep_+(H(\Lambda))\) as full subcategories. By using the same contact isotopy \(H\), we have a functor \(\bar{F}:\mathcal{V}(V,\mathfrak{f}'\cup H(\Lambda^+))\to\mathcal{V}(V,\mathfrak{f}\cup \Lambda^+)\), whose essential image contains \(\mathcal{R}ep_+(\Lambda)\). We also have that \(\bar{F}\) sends linking spheres of \(H(\Lambda^+)\) to linking spheres of \(\Lambda^+\) (see Section \ref{section wrapping exact triangle}). Since removing the stop \(\Lambda^+\) amounts to quotienting by the linking spheres of \(\Lambda^+\) (Theorem \ref{stop removal}), \(\bar{F}\) descends to the functor \(F:\mathcal{V}(V,\mathfrak{f}')\to\mathcal{V}(V,\mathfrak{f})\). Since the essential image of \(\bar{F}\) contains \(\mathcal{R}ep_+(\Lambda)\) (which contains \(\Lambda\) as an \(R\)-brane), we can deduce that the essential image of \(F\) contains \(\Lambda\). Repeating this argument for all \(\Lambda\) proves that \(F\) is essentially surjective.
\end{proof}

\subsection{Wrapping exact triangle}\label{section wrapping exact triangle}
Partially wrapped bilinearized contact cohomology is invariant under Legendrian isotopies away from the stop. We would like to quantify what happens when the isotopy does go through the stop. This is too complicated in general. However, when the stop is itself (locally) a Legendrian submanifold, we have a wrapping exact triangle.

Let \(p\in \mathfrak{f}\) be a point around which locally \(\mathfrak{f}\) is a Legendrian submanifold of \(V\). Let \(\Lambda\subset V\setminus \mathfrak{f}\) be a Legendrian \(R\)-brane, and let \(\Lambda^w\subset V\setminus \mathfrak{f}\) be a Legendrian \(R\)-brane such that there is a positive isotopy from \(\Lambda\) to \(\Lambda\) in \(V\setminus(\mathfrak{f}\setminus p)\) which goes through \(p\) once.

\begin{figure}
   \centering
	\includegraphics[scale=0.6]{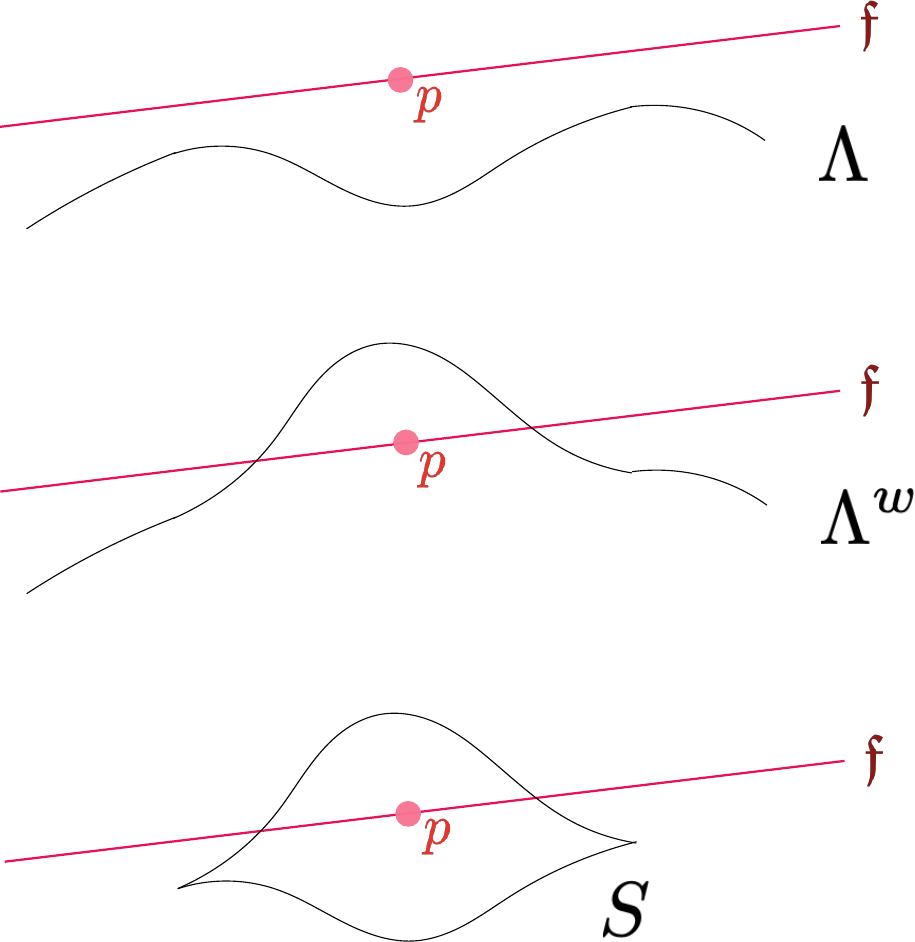}
	\caption{A 1-dimensional example of \(\Lambda,\Lambda^w,S\).}
	\label{figure linking sphere}
\end{figure}

Assume \(\dim V\geq 3\). Let \(S\) be a Legendrian linking sphere of \(\mathfrak{f}\) around the point \(p\), whose Chekanov-Eliashberg DGA is \(\F[t,t^{-1},b]\) if \(\dim(S)=1\), and is \(\F[b]\) if \(\dim(S)>1\). We equipped \(S\) with the unique \(\dim(\epsilon_{\Lambda})\)-dimensional representation \(i:\mathcal{A}(S)\to V^{\Lambda}\) (same vector space as the representation \(\epsilon_{\Lambda}\)), where \(i(b)=0\), and if \(\dim(S)=1\), \(i(t)=1\) or \(i(t)=-1\) (whether \(i(t)=1\) or \(i(t)=-1\) depends on the spin structure we pick on the \(S^1\)). See Figure \ref{figure linking sphere} for the front projection of \(\mathfrak{f},\Lambda,\Lambda',S\).

There is the case where \(V=J^1 *=\R\) is the 1-jet space of a point. In this case, the linking sphere is 2 points, with the Maslov potential of the top point being 1 more than the Maslov potential of the lower point. Here \(\mathcal{A}(S)=\F[b]\), and we set \(i(b)=\id\).

\begin{thm}\label{wrapping exact triangle}
	(Wrapping exact triangle) There is an exact triangle
	\[
	\Lambda\rightarrow\Lambda^w\rightarrow S
	\]
	in the category \(\mathcal{V}(V,\mathfrak{f})\). The map \(\Lambda\rightarrow\Lambda^w\) is the continuation map.
\end{thm}

\begin{figure}
   \centering
	\includegraphics[scale=0.6]{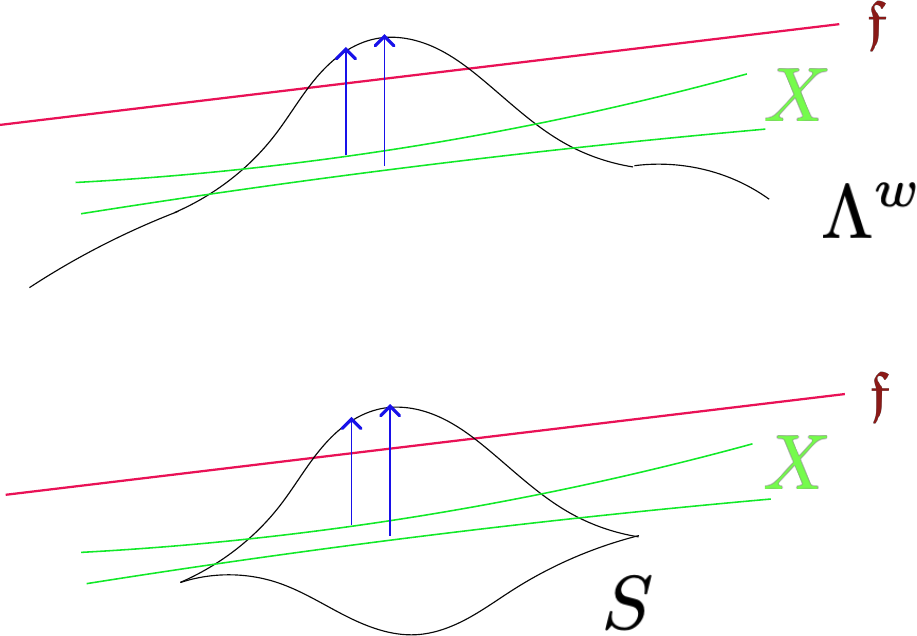}
	\caption{Reeb chords near \(p\) from \(S\) to \(\Lambda^w\) (top) and from \(X\) to \(S\) (bottom) are in bijection.}
	\label{figure X to S chords}
\end{figure}

\begin{proof}
	For the special case where \(\dim V= 1\), the theorem holds trivially.
	
	Now consider the more general case where \(\dim V\geq3\). To prove the theorem, it suffices to prove the statement in the category \(\mathcal{O}\) before localization, since localization preserves exact triangles.

	Let \(X\) be a Legendrian \(R\)-brane that does not intersect \(\mathfrak{f}\). To prove an exact triangle in \(\mathcal{V}(V,\mathfrak{f})\), it suffices to assume that in the poset \(\mathcal{O}\), \(X<\Lambda<\Lambda^w<S\), and prove the existance of an exact triangle
	\[
	LCH(X,\Lambda)\rightarrow LCH(X,\Lambda^w)\rightarrow LCH(X,S).
	\]

	Since this exact triangle is a statement in the cohomology category, we only need to prove it up to a contact isotopy. So we can arrange for the following:
	\begin{enumerate}
		\item \(\Lambda\) and \(\Lambda^w\) are \(C^2\)-close, and agree outside of a Darboux neighborhood \(W\) of \(p\), and that \(S\) is contained in this Darboux neighborhood.
		\item For a neighborhood \(U\) of \(p\), \(X\cap U\) is a disjoint union \(U\cap \bigcup_{i} \Gamma_{c_i}\) of 1-jet graphs of constant functions on \(M\).
		\item \(U\subset W\), and \(S\) is disjoint from \(U\).
		\item The set of Reeb chords from \(X\) to \(\Lambda^w\) is in bijection with union of the set of Reeb chords from \(X\) to \(\Lambda\) and from \(X\) to \(S\). See Figure \ref{figure X to S chords}.
		\item \(W\) is small enough, so that the Reeb chords from \(X\) to \(\Lambda^w\) that are contained in \(W\) (which are the ones that correspond to the Reeb chords from \(X\) to \(S\)) are the shortest Reeb chords from \(X\) to \(\Lambda^w\). By an energy argument, this assumption implies that holomorphic disks on \(\Lambda^w\cup X\) with a positive puncture at a Reeb chord contained in \(W\) are entirely contained within \(W\). Furthermore, these disks are in bijection with the holomorphic disks on \(S\cup X\).
	\end{enumerate}

	With these assumptions, we see that there is a short exact sequence of chain complexes 
	\[
	0\rightarrow LCC(X,\Lambda)\rightarrow LCC(X,\Lambda^w)\rightarrow LCC(X,S)\rightarrow 0,
	\]
	where the arrow \(LCC(X,\Lambda)\rightarrow LCC(X,\Lambda^w)\) is given by inclusion of Reeb chords, which, on the cohomology level, is equal to the map given by composition with the continuation element. The theorem then follows by the Yoneda lemma.
\end{proof}

The wrapping exact triangle has the following immediate corollary:

\begin{thm}\label{stop removal}
	(Stop removal) Let \(\mathfrak{f}'\subset\mathfrak{f}\) be a closed subset so that \(\mathfrak{f}\setminus \mathfrak{f}'\) is an isotropic submanifold. Then,
	\[
	\mathcal{V}(V,\mathfrak{f}')\cong \mathcal{V}(V,\mathfrak{f})/\mathcal{D},
	\]
	where \(\mathcal{D}\) denotes the set of linking spheres of Legendrian points of \(\mathfrak{f}\setminus \mathfrak{f}'\).
\end{thm}

\subsection{Empty stop}

We study the structure of the wrapped contact Fukaya category of a 1-jet space, which is simply the partially wrapped contact Fukaya category with empty stop.

\begin{thm}\label{compact wcfuk with empty stop}
	Let \(M\) be a closed, orientable smooth manifold. Then, \(\mathcal{V}(J^1M,\emptyset)\) is equivalent to the derived category of finite dimensional local systems on \(M\), which is generated by Legendrian \(R\)-branes whose underlying Legendrian is a section of \(J^1M\).
\end{thm}

\begin{proof}
	When the stop is empty, all wrapping can be done using the vertical Reeb field \(\frac{\partial}{\partial z}\). Then, we see that in \(\mathcal{V}(J^1M,\emptyset)\), when computing \(m_k(\Lambda_0,\ldots,\Lambda_k)\), the wrapping needed is so that in the \(z\) coordinate, \(\Lambda_0\) is entirely below \(\Lambda_1\), which is entirely below \(\Lambda_2\), and so on. In this case, instead of counting Reeb chords and holomorphic disks in the symplectization, we can equivalently count intersection points and holomorphic disks in the Lagrangian projection in \(T^*M\). This realized any object of \(\mathcal{V}(J^1M,\emptyset)\) as a compact immersed Lagrangian in \(T^*M\) equipped with a compact local system and a bounding cochain. Any compact immersed Lagrangian in \(T^*M\) equipped with a compact local system and a bounding cochain can be realized as a finite dimensional module over the wrapped Fukaya category of \(T^*M\). This can be further realized as a fully faithful functor \(\mathcal{V}(J^1M,\emptyset)\to\Mod_{fd}\mathcal{W}(T^*M)\) by \cite{CDGG2} (arguments in \cite{CDGG2} essentially prove that the Morita equivalence class of \(\mathcal{W}(T^*M)\) does not change if in the definition of \(\mathcal{W}(T^*M)\), we add immersed Lagrangians equipped with representations of the Chekanov-Eliashberg DGAs of their Legendrian lifts. The only gap is that they did not carefully provide Hamiltonian perturbations needed for defining the \(A_\infty\)-structure of the enlarged category, but this is done using localization in the present paper). By Abouzaid's generation theorem \cites{Abouzaid,Abouzaid2}, \(\Mod_{fd}\mathcal{W}(T^*M)=\Mod_{fd}C_{-*}(\Omega M)\), which is the derived category of finite dimensional local systems on \(M\). By \cite{Abouzaid3}*{Lemma A.1, A.2, A.3}, \(\Mod_{fd}\mathcal{W}(T^*M)\) is generated by the sections \(M\subset T^*M\) equipped with all finite dimensional local systems. Since the 0-sections \(M\subset J^1M\) equipped with all finite dimensional representations are mapped to sections \(M\subset T^*M\) equipped with all finite dimensional local systems, the functor \(\mathcal{V}(J^1M,\emptyset)\to\Mod_{fd}\mathcal{W}(T^*M)\) is essentially surjective.
\end{proof}

\begin{thm}\label{noncompact wcfuk with empty stop}
	Let \(M\) be a noncompact smooth manifold. Then, \(\mathcal{V}(J^1M,\emptyset)\) is trivial.
\end{thm}

\begin{proof}
	Since \(M\) is a noncompact manifold, it admits a smooth function \(f:M\to\R\) with no critical point. Since \(f\) has no critical point, the contact Hamiltonian isotopy on \(J^1M\) given by \((f,df)\) has the property that for any compact subset \(\Lambda\subset J^1M\), there exists some \(t\in\R\) such that the time \(t\) flow of \(\Lambda\) has disjoint Lagrangian projection with \(\Lambda\). This shows that any object in \(\mathcal{V}(J^1M,\emptyset)\) has trivial self-hom, which means that \(\mathcal{V}(J^1M,\emptyset)\) is trivial.
\end{proof}

\section{Partially wrapped microlocal sheaf categories}\label{section sheaf categories}

For the basics of microlocal sheaf theory, we refer the reader to \cite{KS}. In this paper, we use Kuo's wrapped sheaves model \cite{Kuo}. All sheaves in this paper take values in the derived category of modules over a field \(\F\).

We write \[S^*(M\times\R):=T^*(M\times\R)/\R^+\] for the unit cotangent bundle of \(M\times\R\). Decompose \[S^*(M\times \R)=S^+(M\times \R)\cup S^-(M\times \R),\] where 
\[
S^+(M\times \R):=\{\xi\in T^*(M\times\R):\xi(0,\frac{\partial}{\partial z})\geq 0\}/\R^+,
\]
and
\[S^-(M\times \R):=\{\xi\in T^*(M\times\R):\xi(0,\frac{\partial}{\partial z})\leq 0\}/\R^+.
\]
There is a contactomorphism
\[
J^1M\cong S^*(M\times\R)\setminus S^-(M\times\R)\subset S^+(M\times\R)
\]
which is an identification we use implicitly for the rest of this paper.

\begin{defn}
	Let \(M\) be a real analytic manifold. A \emph{Legendrian \(S\)-brane} (\(S\) is for sheaf) consists of the data of a compact Legendrian \(\Lambda\subset J^1M\subset S^*(M\times\R)\), and a complex of sheaves on \(M\times\R\) with perfect stalks whose microsupport is contained in \(\Lambda\). We require that the complex of sheaves has vanishing cohomology near \(M\times\{\infty\}\). We will often abuse notation and simply denote a Legendrian \(S\)-brane by its underlying Legendrian \(\Lambda\).

	For the notion of \emph{isotopies of Legendrian \(S\)-branes}, see \cite{Kuo}*{Section 3}. Note that Kuo simply calls them sheaves.

	The morphism space (in the unwrapped derived category of sheaves) between Legendrian \(S\)-branes \(\Hom_S(\Lambda,\Lambda')\) is simply the morphism space in the derived category of sheaves.
\end{defn}

\begin{defn}
	Let \(M\) be a real analytic manifold, and let \(\mathfrak{f}\subset J^1M\) be a closed subset. Consider the category \(\operatorname{\mathfrak{w}sh}_{\mathfrak{f}\cup S^-(M\times\R)}(M\times \R)\) of wrapped sheaves as in \cite{Kuo}*{Definition 4.1}. In this paper, since we do not consider other sheaf categories, we let \(\Sh(J^1M,\mathfrak{f})\) denote the full subcategory of \(\operatorname{\mathfrak{w}sh}_{\mathfrak{f}\cup S^-(M\times\R)}(M\times \R)\) consisting of Legendrian \(S\)-branes in \(J^1M\setminus \mathfrak{f}\), which we call the \emph{partially wrapped microlocal sheaf category}.
\end{defn}

\begin{rem}
	\(\operatorname{\mathfrak{w}sh}_{\mathfrak{f}\cup S^-(M\times\R)}(M\times \R)\) has more objects than \(\Sh(J^1M,\mathfrak{f})\), since it also contains sheaves whose stalks are not perfect, and sheaves that do not vanish near \(M\times\{\infty\}\).
\end{rem}

Note that similar to the partially wrapped contact Fukaya category, we have that in the partially wrapped microlocal sheaf category, \[H^*\Sh(\Lambda_1,\Lambda_2)=\varinjlim_i\Hom_S(\Lambda_1,\Lambda_2^t),\] where \(\Lambda_2^t\) is a positive cofinal wrapping of \(\Lambda_2\).

Observe that given an unknot \(S\subset J^1M\) with standard front projection, and a positive integer \(n\), there is a unique sheaf of microlocal rank \(n\) with microsupport on \(S\), up to grading shift. We can therefore view a linkng sphere at a Legendrian point of the stop as an \(S\)-brane.

Similar to the partially wrapped contact Fukaya categories, there is a wrapping exact triangle in the sheaf category \(\Sh(J^1M,\mathfrak{f})\). Let \(p\in \mathfrak{f}\) be a point around which locally \(\mathfrak{f}\) is a Legendrian submanifold of \(V\). Let \(\Lambda\subset V\setminus \mathfrak{f}\) be a Legendrian \(S\)-brane, and let \(\Lambda^w\subset V\setminus \mathfrak{f}\) be a Legendrian \(R\)-brane such that there is a positive isotopy from \(\Lambda\) to \(\Lambda\) in \(V\setminus(\mathfrak{f}\setminus p)\) which goes through \(p\) once.

\begin{thm}[\cite{Kuo}*{Definition 4.15, Remark 4.16}]\label{sheaf wrapping exact triangle}
	(Wrapping exact triangle) There is an exact triangle
	\[
	\Lambda\rightarrow\Lambda^w\rightarrow S
	\]
	in the category \(\Sh(V,\mathfrak{f})\). The map \(\Lambda\rightarrow\Lambda^w\) is the continuation map.
\end{thm}

\begin{thm}\label{sheaf with empty stop}
	If \(M\) is compact, \(\Sh(J^1M,\emptyset)\) is equivalent to the derived category of finite dimensional local systems on \(M\). If \(M\) is noncompact, \(\Sh(J^1M,\emptyset)\) is trivial.
\end{thm}

\begin{proof}
	If \(M\) is compact, then we see that after completely wrapping a Legendrian with a finite rank sheaf, the result is a finite rank sheaf on \(M\times\R\), which is equivalently a finite rank sheaf on \(M\). By taking the Legendrian to be a section, we see that we can obtain any finite rank sheaf on \(M\).

	If \(M\) is noncompact, then any sheaf on a compact Legendrian that vanishes at positive infinity automatically vanishes at infinity. Therefore, after wrapping, the sheaf vanishes everywhere.
\end{proof}

\section{Microlocal Morse theory of partially wrapped contact Fukaya categories}\label{section equivalence}

In this section, we prove the main theorem of this paper, Theorem \ref{main theorem}. We prove it by computing both categories explicitly. At a first glance, explicit computation might seem intractible, since cofinal wrapping sequences are in general quite wild. In general, it is not even possible to find generators of the partially wrapped contact Fukaya category with easy to understand cofinal wrapping sequences. However, it turns out that if we first refine the stop (Section \ref{section stop refinement}), then we can find generators with easy to describe cofinal wrapping sequences. If we pick the stop refinement carefully enough, we can prove that the subcategory of generators is formal (Section \ref{section morphisms and compositions}). In this strategy, one caveat is that the generators coming from Theorem \ref{compact wcfuk with empty stop} are local systems, between which the restriction of the partially wrapped contact Fukaya category is not formal. We deal with these section generators by first enlarging the category so that we can use a different set of generators between which the category is formal, and then characterize the local systems by the Yoneda lemma.

\subsection{Stop refinement}\label{section stop refinement}

Let \(M\) be a real analytic manifold. A closed subanalytic set \(\mathfrak{f}\subset J^1M\) is called isotropic if for some, or equivalently for every, cover of \(\mathfrak{f}\) by locally closed \(C^1\) manifolds, all of them are isotropic.

Assume that \(\mathfrak{f}\subset J^1M\) is closed, subanalytic, isotropic, and bounded in the vertical direction. Without loss of generality, assume \(\mathfrak{f}\subset T^*M\times[0,1]\). Let \(\mathcal{T}\) be a subanalytic stratification of \(M\times[0,1]\) so that the front projection \(\Pi(\mathfrak{f})\) is a union of \(\mathcal{T}\)-strata. We may assume that no strata of \(\mathcal{T}\) has any vertical tangency, by refining the strata.

Write \(p:M\times\R\to M\) for the base projection map. Let \(\mathcal{T}^p\) be a subanalytic stratification of \(X\) such that for each strata \(t\in\mathcal{T}\) that is not top-dimensional, \(p(t)\) is a union of strata in \(\mathcal{T}^p\). By \cite{Cza}, \(\mathcal{T}^p\) can be refined to be a Whitney \(C^d\)-triangulation \((K,H)\) where \(K\) is an abstract simplicial complex, and \(H\) is a homeomorphism \(|K|\to M\) such that \(H\) is a \(C^d\)-embedding on each simplex.

\begin{figure}
   \centering
	\includegraphics[scale=0.8]{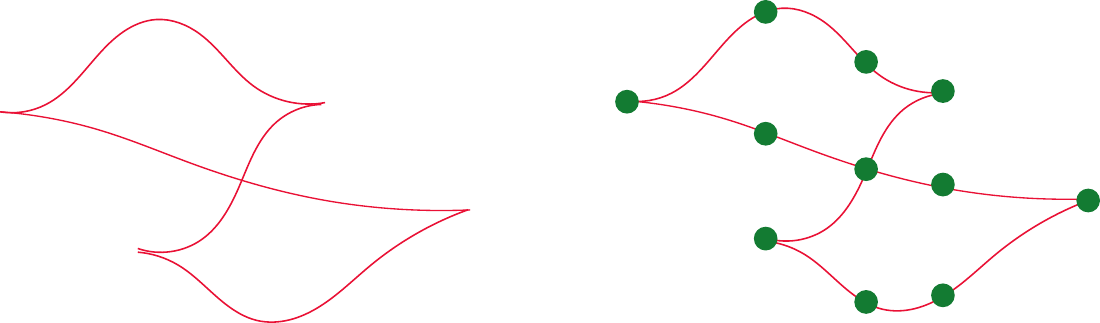}
	\caption{An example of a Legendrian stop \(\mathfrak{f}\) (left) and its refinement \(\mathfrak{S}\) (ight).}
	\label{figure f to S example}
\end{figure}

Let \(\mathcal{S}\) be the subanalytic Whitney stratification of \(M\times[0,1]\) defined as follows: each \(t\in\mathcal{T},t'\in\mathcal{T}^p\) give rise to a stratum \(t\cap (t'\times\R)\in\mathcal{S}\). Write \(\mathfrak{S}:=\bigcup_{s\in\mathcal{S}} N^* s\cap J^1M\), and \(\mathfrak{S}_s:=\bigcup_{t\subset \bar{s}}N^*t\cap J^1M\). Note that for any top dimensional strata \(s\), \(N^*s\) is empty. See Figure \ref{figure f to S example} for the front projection of a stop refinement from \(\mathfrak{f}\) to \(\mathfrak{S}\) in an example in \(J^1\R\).

For the purpose of simplifying arguments in the case where \(M\) is compact, we will also make use of the stop \(\mathfrak{g}:=\bigcup_{\Delta\in K}N^*(\Delta\times \{2\})\cap J^1M\). Write \(\mathfrak{g}_\Delta:=\bigcup_{\Delta'\subset \bar{\Delta}}N^*(\Delta'\times\{2\})\cap J^1M\).

\begin{figure}
   \centering
	\includegraphics[scale=0.8]{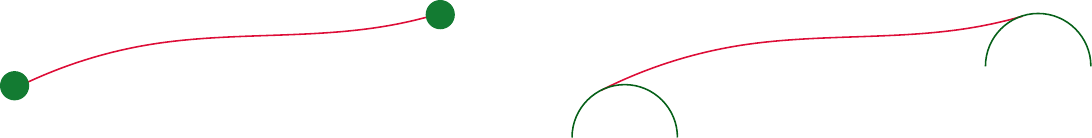}
	\caption{A 1-dimensional example of deforming \(\mathfrak{S}_s\) (left) to \(\mathfrak{S}_s^+\) (right).}
	\label{figure 1 dim deformation}
\end{figure}

\begin{figure}
   \centering
	\includegraphics[scale=0.8]{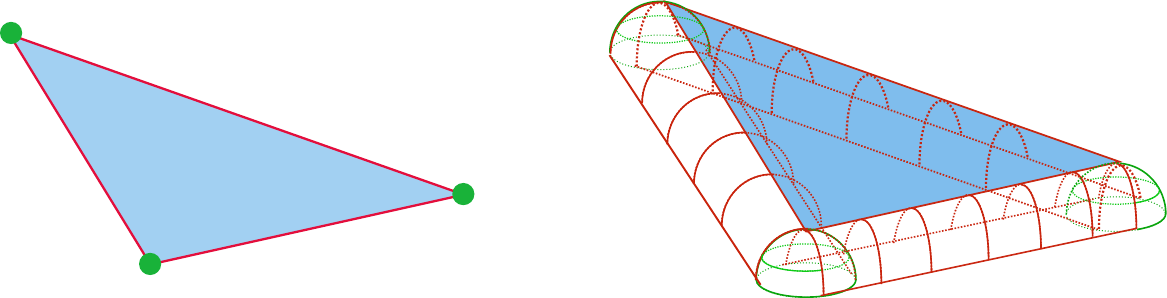}
	\caption{A 2-dimensional example of deforming \(\mathfrak{S}_s\) (left) to \(\mathfrak{S}_s^+\) (right).}
	\label{figure 2 dim deformation}
\end{figure}

Equipping \(M\) with a Riemannian metric, and \(\R\) with the standard Riemannian metric. This gives a product Riemannian metric on \(M\times\R\), which defines a Reeb field on \(S^*(M\times\R)\). For any subset \(X\subset J^1M\), write \(X^+\) for a small pushoff of \(X\) along the Reeb field given by the metric. We have that \(p(\mathfrak{S}_s^+)\) is a regular neighborhood of \(p(\mathfrak{S}_s)\), and is homeomorphic to \(B^n\). See Figure \ref{figure 1 dim deformation} for the front projection of an example of \(\mathfrak{S}_s\) deformed into \(\mathfrak{S}_s^+\) when \(\dim M=1\). See Figure \ref{figure 2 dim deformation} for a 2-dimensional example.

\begin{assumption}\label{assumption}
We make the assumption that \(p(\mathfrak{S}_s^+)\) is (after smoothing corners) diffeomorphic to \(B^n\). More precisely, we assume that for the stop \(\mathfrak{f}\) that we started with, a choice that results in all \(p(\mathfrak{S}_s^+)\) being diffeomorphic to \(B^n\) exists. This is only an assumption when \(n=4,5\), since by the \(h\)-cobordism theorem, any smooth manifold homeomorphic to \(B^n\) for \(n\geq 6\) is diffeomorphic to \(B^n\), and the same is true if \(n\leq 3\). This assumption makes it easy to define \(T_s\) and \(Q_\Delta\). Note that this assumption in dimension \(4\) follows from the 4-dimensional smooth Schoenflies conjecture, and in dimension \(5\) follows from the 4-dimensional smooth Poincaré conjecture.
\end{assumption}

\begin{rem}
	The assumption that \(p(\mathfrak{S}_s^+)\) is (after smoothing corners) diffeomorphic to a ball is likely not necessary. It is possible that in our setting, they are diffeomorphic to \(B^n\) without any assumption, but in general, exotic smooth structures on balls of dimensions \(4,5\) exist iff exotic smooth structures on \(S^4\) exist, which is question posed by the 4-dimensional smooth Poincaré conjecture. In the case of \(n=4\), the ball can be embedded in standard \(\R^4\), so the 4-dimensional smooth Schoenflies conjecture suffices. Even if \(p(\mathfrak{S}_s^+)\) in our setting is not diffeomorphic to a ball, it is still likely that there exists a stabilization trick like the one in \cite{GPS3} that allows us to conclude all the same theorems later in this paper.
\end{rem}

Let \(s\in\mathcal{S}\) be a strata that is not top dimensional. We define the subset \(B_s\subset \mathfrak{S}_s^+\) as the conormal to the top most layer in the front projection, that is, in the front projection \(\Pi(\mathfrak{S}_s^+)\), given a point \(x\in M\), if the intersection \(\Pi(\mathfrak{S}_s^+)\cap (x\times\R)\) has more than 1 point, then we only keep the point that is the highest.

In the following discussion, in order to declutter the notation, we identify a smooth function from an embedded ball of \(M\) to \(\R\) with both its graph in \(M\times\R\) and its \(1\)-jet graph. Let \(S_{s}\) be a linking sphere of a point on \((N^* s\cap J^1M)^+\). A cofinal wrapping of \(S_{s}\) in \(J^1M\setminus (\mathfrak{S}_s\cup \mathfrak{g}_{p(s)})^+\) can be given as follows. Let \(D_{s},D_{s}'\) denote the lower and upper sheets of a standard 2-sheet front projection of the Legendrian unknot \(S_{s}\). Let \(D^n\) be the standard closed disk of dimension \(n\). Choose a positive isotopy \(\phi_t:D^n \hookrightarrow J^1M\setminus (\mathfrak{S}_s\cup \mathfrak{g}_{p(s)})^+\) so that \(\phi_0(D^n)=D_s\) and \(\lim_{t\to \infty}\phi_t(D^n)=B_s\) (we can pick the isotopy so that it never intersects with \(\mathfrak{S}_s^+\setminus B_s\)). Pick a smooth family of increasing Morse function on \(f_t:p(\phi_t(D^n))\to\R\) such that: the derivatives of \(f_t\) of all orders vanish on the boundary, all \(f_t\) take positive values on the interior, and \(f_0=D'_s\) and \(f_t\) converges pointwise to \(2-B_s\) as functions. Let \(\psi_t:D^n\hookrightarrow M\times\R\) be an isotopy defined by \(\psi(D^n)=\phi_t(D^n)+f_t\). Then, \(\phi\) and \(\psi\) agree on their boundary, so they together define an isotopy of the Legendrian unknot \(S_{s}\). We can choose \(f_t\) so that the isotopy is in \(J^1M\setminus (\mathfrak{S}_s\cup \mathfrak{g}_{p(s)})^+\). Since this positive isotopy eventually escapes all compact subsets of \(J^1M\setminus (\mathfrak{S}_s\cup \mathfrak{g}_{p(s)})^+\), it is a cofinal positive wrapping.

Since we assumed \(\mathfrak{f}\) is bounded in the vertical direction, for each stratum \(s\in\mathcal{S}\), there are only finitely many strata of \(\mathcal{S}\) above a neighborhood of \(s\). So, the cofinal wrapping of \(S_{s}\) defined above eventually goes through the conormals of all those finitely many strata, and the Legendrian isotopy class of \((S_{s})_t\) in \(J^1M\setminus (\mathfrak{S}\cup \mathfrak{g})^+\) eventually stabilizes.

\begin{figure}
   \centering
	\includegraphics[scale=0.5]{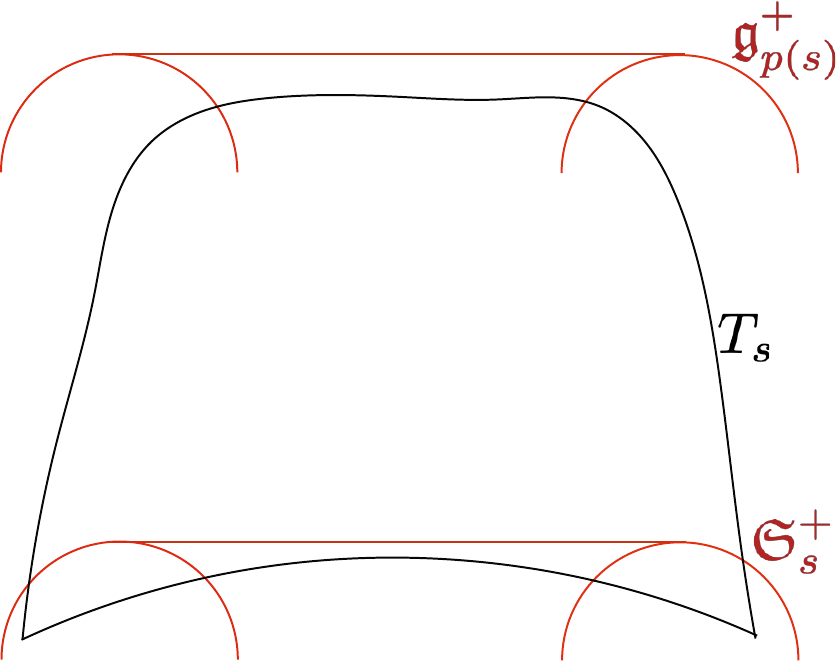}
	\caption{A 1-dimensional example of \(T_s\).}
	\label{figure Ts big}
\end{figure}

\begin{figure}
   \centering
	\includegraphics[scale=0.5]{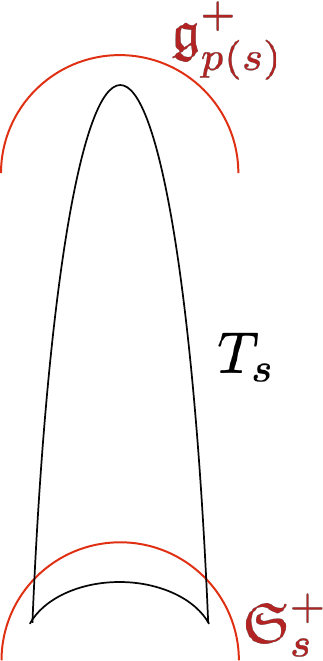}
	\caption{Another 1-dimensional example of \(T_s\).}
	\label{figure Ts small}
\end{figure}

\begin{defn}
	We pick a representative of this Legendrian isotopy class, and call it \(T_{s}\). For ease of notation, if \(s\) is top-dimensional, then we set \(T_s:=\emptyset\). See Figures \ref{figure Ts big} and \ref{figure Ts small} for \(1\)-dimensional examples of \(T_s\).
\end{defn}

It is easy to find a cofinal wrapping of \(T_{s}\) in \(J^1M\setminus (\mathfrak{S}\cup\mathfrak{g})^+\): we simply continue the isotopy we used to obtain \(T_{s}\).

For a simplex \(\Delta\in K\), we define \(Q_{\Delta}\) in a similar way to the above discussion: Let \(S_{\Delta}\) be a linking sphere of a point on \((N^* (\Delta\times\{2\})\cap J^1M)^+\). A cofinal wrapping of \(S_{\Delta}\) in \(J^1M\setminus (\mathfrak{g}_{\Delta})^+\) can be given as follows. Let \(D_{\Delta},D_{\Delta}'\) denote the lower and upper sheets of a standard 2-sheet front projection of the Legendrian unknot \(S_{\Delta}\). Let \(D^n\) be the standard closed disk of dimension \(n\). Choose a positive isotopy \(\phi_t:D^n \hookrightarrow J^1M\setminus (\mathfrak{g}_{\Delta})^+\) so that \(\phi_0(D^n)=D_\Delta\) and \(\lim_{t\to \infty}\phi_t(D^n)=B_{\Delta\times\{2\}}\). Pick a Morse function \(f:D^n\to\R\) such that the derivatives of \(f\) of all orders vanish on the boundary, and \(f\) takes positive values on the interior. Let \(\psi_t:D^n\hookrightarrow M\times\R\) be an isotopy defined by \(\psi(D^n)=\phi_t(D^n)+tf\). Then, \(\phi\) and \(\psi\) agree on their boundary, so they together define an isotopy of the Legendrian unknot \(S_{s}\). It is not hard to see that we can choose \(f\) so that the isotopy is in \(J^1M\setminus (\mathfrak{g}_{\Delta})^+\). Since this isotopy eventually escapes all compact subsets of \(J^1M\setminus (\mathfrak{g}_{\Delta})^+\), it is a cofinal wrapping. In fact this cofinal wrapping is also contained in \(J^1M\setminus (\mathfrak{S}\cup\mathfrak{g})^+\). The Legendrian isotopy class of this cofinal wrapping eventually stablizes.

\begin{figure}
   \centering
	\includegraphics[scale=0.5]{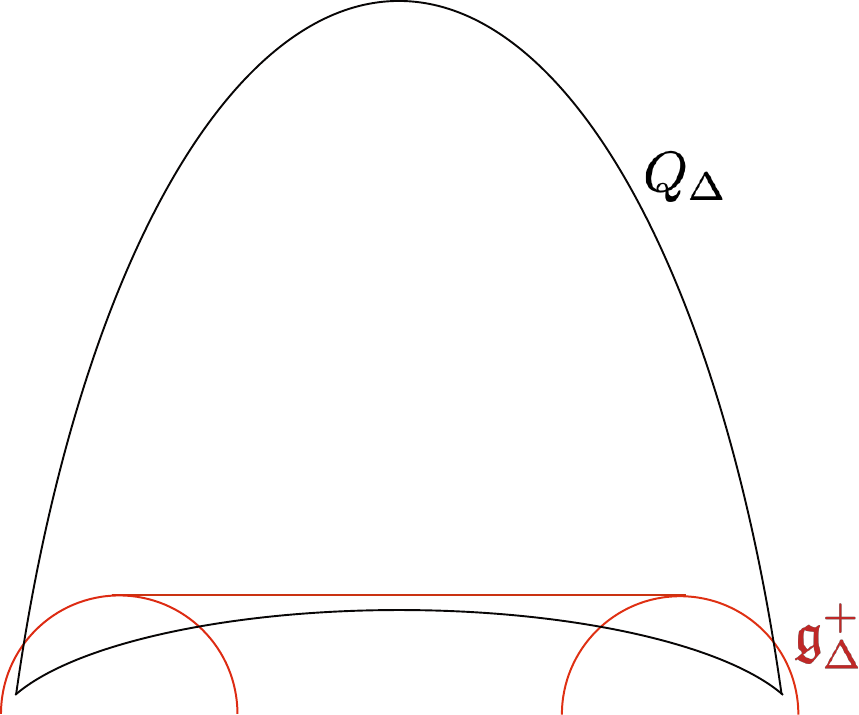}
	\caption{A 1-dimensional example of \(Q_\Delta\).}
	\label{figure Q_delta}
\end{figure}

\begin{figure}
   \centering
	\includegraphics[scale=0.5]{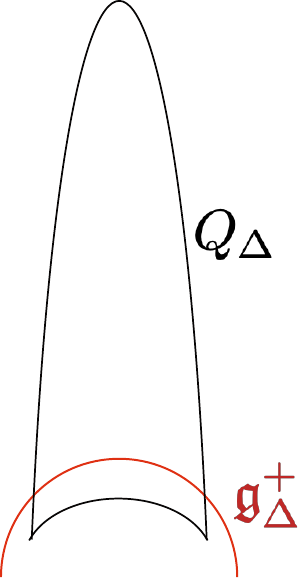}
	\caption{Another 1-dimensional example of \(Q_\Delta\).}
	\label{figure Q_Delta small}
\end{figure}

\begin{defn}
	We pick a representative of this Legendrian isotopy class, and call it \(Q_\Delta\). See Figures \ref{figure Q_delta} and \ref{figure Q_Delta small} for \(1\)-dimensional examples of \(Q_\Delta\).
\end{defn}

\subsection{Generators}

When dealing with the partially wrapped contact Fukaya category, we equip \(T_{s}\)'s, \(Q_\Delta\)'s, and linking spheres (that have a standard front projection consisting of 2 sheets and a sphere worth of cusps) with Legendrian \(R\)-brane structures by equipping them with the unique 1-dimensional representation. We equip them with basepoints in the center of their bottom sheets, and arbitrary endpoint paths. In the case where \(\dim(M)=1\), we equip the circles with the null-cobordant spin structure. We equip the bottom sheets of \(T_{s}\)'s and linking spheres of \(\mathfrak{S}\) with Maslov potential \(0\). We equip the bottom sheets of \(Q_{\Delta}\)'s and linking spheres of \(\mathfrak{g}\) with Maslov potential \(1\).

When dealing with the partially wrapped sheaf category, we equip the above Legendrians with Legendrian \(S\)-brane structures given by the unique \(1\)-dimensional sheaf. For \(T_{s}\)'s and linking spheres of \(\mathfrak{S}\), we let the sheaf be concentrated in degree \(0\). For \(Q_{\Delta}\)'s and linking spheres of \(\mathfrak{g}\), we let the sheaf be concentrated in degree \(1\).

All arguments in this subsection use only the wrapping exact triangle, and generation of the category with empty stop by sections equipped with finite dimensional local systems. These two ingredients are true in both the partially wrapped contact Fukaya category, and the partially wrapped microlocal sheaf category. Therefore, the results of this subsection are true in both categories. Thus we write \(\mathcal{C}\) for either \(\mathcal{V}\) or \(\Sh\). We write \emph{brane} for either an \(R\)-brane or an \(S\)-brane. 

\begin{rem}
	In the proposition below, note that although \(T_{s}\) is well-defined, \(S_{s}\) is not well-defined as a Legendrian in \(J^1M\setminus\mathfrak{S}\), since linking spheres of different points on \(N^*s\cap J^1M\) might not be isotopic. An isotopy from one such linking sphere to another might need to pass through the conormal of strata in \(\Star(s)\). Depending on context, when we write \(S_{s}\), we might mean either the collection of all isotopy classes of linking spheres on \(N^*s\cap J^1M\), or the linking sphere of any point on \(N^*s\cap J^1M\), or the linking sphere of some particular point on \(N^*s\cap J^1M\).
\end{rem}

\begin{prop}\label{prop t=s}
	The subcategory of \(\mathcal{C}(J^1M,\mathfrak{S})\) generated by all linking spheres on \(\mathfrak{S}\) is equal to the subcategory generated by \(\{T_{s}:s\in\mathcal{S}\}\).
\end{prop}

\begin{proof}
	Let \(0\leq i\leq n\). For any \(i\)-dimensional stratum \(s\), the process of isotoping any \(S_{s}\rightsquigarrow T_{s}\) passes through finitely many strata of dimension in \(\{i+1,\ldots,n\}\), and the \(i\)-dimensional strata that are directly above \(s\).

	To prove the proposition, it suffices to prove that any linking sphere \(S_{s}\) is a twisted complex of several \(T_{t}\)'s. To prove this, we proceed by two induction arguments. For any \(n\)-dimensional stratum \(s\subset M\times 1\), we have that \(T_{s}\) is isotopic to \(S_{s}\). Assume that an \(n\)-dimensional stratum \(s\) has \(k\) \(n\)-dimensional strata \(t_1,\ldots,t_k\) directly above it, and suppose the subcategory generated by all linking spheres on \(\bigcup N^* t_{i}\cap J^1M\) is equal to the subcategory generated by \(\{T_{t_1},\ldots,T_{t_k}\}\). Then, by the wrapping exact triangle, any linking sphere \(S_{s}\) on a smooth point of \(N^* s\cap J^1M\) is isomorphic to a twisted complex of \(T_{s},S_{1},\ldots,S_{m}\) (where \(S_1,\ldots,S_m\) is some sequence of linking spheres on \(\bigcup N^* t_{i}\cap J^1M\)), which by the induction hypothesis, is isomorphic to twisted complex of \(T_{s},T_{t_1},\ldots,T_{t_k}\).

	Let \(s\) be an \(i\)-dimensional stratum with \(k\) \(i\)-dimensional strata \(t_1,\ldots,t_k\) directly above it. Assume that the subcategory generated by \(S_{t}\) for all strata of dimension \(i+1,\ldots,n\) is equal to the subcategory generated by \(T_{t}\) for all strata of dimension \(i+1,\ldots,n\). Then, by the wrapping exact triangle, similar to above, any linking sphere \(S_{s}\) is a twisted complex of \(T_{s},T_{t_{1}},\ldots,T_{t_{k}},T_{l_1},\ldots,T_{l_m}\) where \(l_1,\ldots,l_m\) are strata of dimension greater than \(i\). So, by induction, any \(S_{s}\) is a twisted complex of several \(T_{t}\)'s, which finishes the proof.
\end{proof}

\begin{rem}
	The reason we use two induction arguments above instead of one is that given a particular \(\mathcal{S}\), there might be arbitrarily many \(i\)-dimensional strata directly above an \(i\)-dimensional stratum.
\end{rem}

\begin{cor}
	\(\mathcal{C}(J^1M,\mathfrak{S})\) is generated by the union
	\[
	\{T_{s}:s\in \mathcal{S}\}\cup\{\Sigma:\text{the underlying Legendrian of }\Sigma\text{ is a section below }T^*M\times\{0\}\}.
	\]
\end{cor}

\begin{prop}\label{prop t+q=s}
	The subcategory of \(\mathcal{C}(J^1M,(\mathfrak{S}\cup\mathfrak{g})^+)\) generated by all linking spheres on \(\mathfrak{S}\cup\mathfrak{g}\) is equal to the subcategory generated by \(\{T_s:s\in\mathcal{S}\}\cup\{Q_\Delta:\Delta\in K\}\).
\end{prop}
\begin{proof}
	By an analogous argument to that of Proposition \ref{prop t=s}, we see that the subcategory of \(\mathcal{V}(J^1M,(\mathfrak{S}\cup\mathfrak{g})^+)\) generated by linking spheres on \(\mathfrak{g}\) is equal to the subcategory generated by \(\{Q_\Delta:\Delta\in K\}\). Together with Proposition \ref{prop t=s}, this proves the claim.
\end{proof}

\begin{defn}
	Let \(\mathcal{P}\subset \mathcal{C}(J^1M,(\mathfrak{S}\cup\mathfrak{g})^+)\) be the subcategory whose objects are Legendrians which are above \(T^* M\times\{2\}\).
\end{defn}

\begin{prop}
	Consider the quotient category \(\mathcal{C}(J^1M,(\mathfrak{S}\cup\mathfrak{g})^+)/\mathcal{P}\). Its full subcategory consisting of Legendrians which are below \(T^* M\times\{2\}\) is equivalent to \(\mathcal{C}(J^1M,\mathfrak{S}^+)\).
\end{prop}

\begin{proof}
	The full subcategory of \(\mathcal{C}(J^1M,(\mathfrak{S}\cup\mathfrak{g})^+)\) consisting of Legendrians which are below \(T^* M\times\{2\}\) is clearly equivalent to \(\mathcal{C}(J^1M,\mathfrak{S}^+)\). For any object \(P\in \mathcal{P}\) and any brane \(\Lambda\) below \(T^* M\times\{2\}\), we have that \(H^*\mathcal{C}(P,\Lambda)=0\). So, for \(\Lambda,\Lambda'\) below \(T^* M\times\{2\}\), we have 
	\[
	H^*(\mathcal{C}(J^1M,(\mathfrak{S}\cup\mathfrak{g})^+))(\Lambda,\Lambda')=H^*(\mathcal{C}(J^1M,(\mathfrak{S}\cup\mathfrak{g})^+)/\mathcal{P})(\Lambda,\Lambda'),
	\]
	which proves the claim.
\end{proof}

\begin{prop}\label{quotient category generators}
	The quotient category \(\mathcal{C}(J^1M,(\mathfrak{S}\cup\mathfrak{g})^+)/\mathcal{P}\) is generated by \(\{T_s:s\in\mathcal{S}\}\cup\{Q_\Delta:\Delta\in K\}\).
\end{prop}
\begin{proof}
	After wrapping any Legendrian through stops till it is above \(T^* M\times\{2\}\), it is the \(0\) object by definition. So, \(\mathcal{C}(J^1M,(\mathfrak{S}\cup\mathfrak{g})^+)/\mathcal{P}\) is generated by linking spheres, and so by Proposition \ref{prop t+q=s}, it is generated by \(\{T_s:s\in\mathcal{S}\}\cup\{Q_\Delta:\Delta\in K\}\).
\end{proof}

\subsection{Morphisms and compositions}\label{section morphisms and compositions}

Next we compute the cohomology of morphisms between the above mentioned Legendrian branes. At this point, this is not too hard, since we have already explicitly described wrappings of generators. Define a poset structure on \(\mathcal{S}\) (without the top dimensional strata) as follows. Each \(i\)-dimensional stratum \(s\) is less than all strata \(t\) such that: in the base projection, \(p(s)\subset p(\bar{t})\), and either \(s\subset \bar{t}\) or \(\bar{t}\) is vertically above \(s\). When we say a subset \(X\) is vertically above another subset \(Y\), we mean that for any point \(y\in Y\), there exists a point \(x\in X\) such that \(p(x)=p(y)\), and the vertical coordinate of \(x\) is greater than or equal to the vertical coordinate of \(y\).

\begin{prop}\label{morphisms between generators}
	Equip \(\Sigma\) with Maslov potential \(0\). Then, in the cohomology category \(H^*\mathcal{C}(J^1M,\mathfrak{S})\), we have the following:
	\begin{enumerate}
		\item If \(s\leq t\), then \(H^*\Hom(T_{s},T_{t})=\F\).
		\item If \(s\nleq t\), then \(H^*\Hom(T_{s},T_{t})=0\).
		\item \(H^*\Hom(\Sigma,T_{s})=0\).
		\item If \(s\subset M\times\{0\}\), then \(H^*\Hom(T_{s},\Sigma)=V^\Sigma\).
		\item If \(s\not\subset M\times\{0\}\), then \(H^*\Hom(T_{s},\Sigma)=0\).
		\item If \(p(s)\subset \bar{\Delta}\), then \(H^*\Hom(Q_{\Delta},T_s)=\F\).
		\item If \(p(s)\not\subset \bar{\Delta}\), then \(H^*\Hom(Q_{\Delta},T_s)=0\).
		\item \(H^*\Hom(T_s,Q_{\Delta})=0\).
		\item If \(\Delta\subset\bar{\Delta}'\), then \(H^*\Hom(Q_\Delta,Q_{\Delta'})=\F\).
		\item If \(\Delta\not\subset\bar{\Delta'}\), then \(H^*\Hom(Q_\Delta,Q_{\Delta'})=0\).
	\end{enumerate}
   In (1), (6), (9), the continuation map is identified with \(1\in\F\).
\end{prop}

\begin{figure}
   \centering
	\includegraphics[scale=0.8]{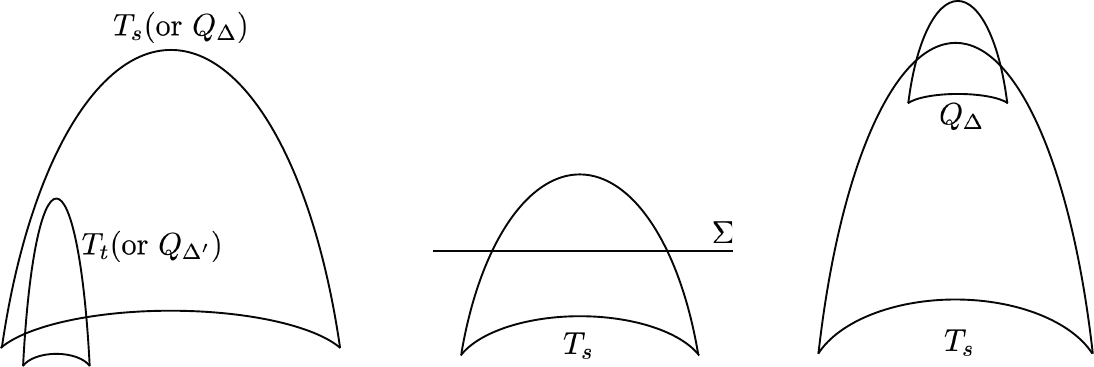}
	\caption{Isotopy types of the union of domain and target after wrapping, used in the computation of morphisms: cases (1), (9) (left), case (4) (middle), case (6) (right).}
	\label{figure morphisms}
\end{figure}

\begin{proof}[Proof for \(R\)-branes]
	We use the fact that in all cases, after suffiently wrapping the second entry, the Legendrian isotopy type of the union of the first and second entry stabilizes. We then use isotopy invariance of morphisms of branes to compute the morphism spaces. See Figure \ref{figure morphisms} for the front projection of \(1\)-dimensional examples of cases (1),(4),(6),(9).

	In (1), (9), after sufficiently wrapping the second entry, the front of the union of the first and second entry is isotopic to the front of 2 parallel copies of the standard unknot. We can arrange so that up to isotopy, there are 3 Reeb chords from an unknot \(U\) to a positive pushoff \(U^+\), with degrees \(0,n,n+1\), these are respectively the short Reeb chord corresponding to an index 0 critical point on \(U\) (which we denote \(\min\)), the short Reeb chord corresponding to an index n critical point on \(U\) (which we denote \(\max\)), and the unique long Reeb chord \(b\). There is a rigid holomorphic strip with positive puncture at \(b\), and negative puncture at \(\max\). There is no other rigid holomorphic curve present, so \(m_1(\max^\vee)=b\). We conclude that \(H^*\Hom(T_{s},T_{t})=\F\), which is generated by the continuation element.

	In (4), (6), after sufficiently wrapping the second entry, we can isotope the union of the two Legendrians so that there is a single Reeb chord from the first entry to the second entry, with Morse index \(0\) and Maslov potential difference \(0\).

	In all other cases, after suffiently wrapping the second entry, either the first or the second entry is a Legendrian unknot that is unlinked from the other entry. So we can isotope it to somewhere far away so that there is no Reeb chord between the two entries.
\end{proof}

\begin{proof}[Proof for \(S\)-branes]
	We use the same fact about the wrappings. In (1), (9), after sufficiently wrapping the second entry, the front of the union of the first and second entry is isotopic to the front of 2 parallel copies of the standard unknot. In this case, the morphism is \(\F\), since the morphisms are sheaves supported on the bounded region by something isotopic to an unknot, and microsupported by something isotopic to an unknot.

	In (4), (6), after suffciently wrapping the second entry, we have that locally around the first entry (which is a standard unknot), the second entry is a horizontal sheet cutting through the top sheet of the first entry with a sheaf supported below the horizontal sheet. The morphism is then a sheaf supported in the intersection of the region below the horizontal sheet and the bounded region in the complement of the front of first entry. The microsupport of the morphism is isotopic to a standard unknot. The result computations in (4) and (6) then follow. Note that in (4), the morphism space \(V^\Sigma\) is naturally identified with the sections of the sheaf of \(\Sigma\) on the interior of the intersection of region below \(\Sigma\) and the bounded region in the complement of the front of \(T_s\).

	In all other cases, after suffiently wrapping the second entry either the first or the second entry is a Legendrian unknot that is unlinked from the other entry. So we can isotope it to somewhere far away so that the sheaves have disjoint support.
\end{proof}

We also need the following computation:

\begin{prop}\label{compositions between generators}
	Let \(X,Y,Z\in \{T_s:s\in\mathcal{S}\}\cup\{Q_\Delta:\Delta\in K\}\). If \(H^*\Hom(X,Y)\) and \( H^*\Hom(Y,Z)\) are both nontrivial, then \[H^*m_2:H^*\Hom(X,Y)\otimes H^*\Hom(Y,Z)\to H^*\Hom(X,Z)\] is given by the multiplication map \(\F\otimes\F\to\F\).

	If \(s,t\subset M\times\{0\}\) and \(s\subset \bar{t}\), then \[H^*m_2:H^*\Hom(T_s,T_t)\otimes H^*\Hom(T_t,\Sigma)\to H^*\Hom(T_s,\Sigma)\] is given by the composition of the multiplication map \(\F\otimes V^\Sigma\to V^\Sigma\) with \(\epsilon_\Sigma(\gamma)\), where \(\gamma\) is the based loop on \(\Sigma\) (here we identify \(\Sigma\) with \(M\)) formed by: the path from \(x\) to the barycentre of \(s\) (from the \(R\)-brane structure of \(T_s\)), concatenated with the path from the barycentre of \(s\) to the barycentre of \(t\), concatenated with the path from the barycentre of \(t\) to \(x\) (from the \(R\)-brane structure of \(T_t\)).
\end{prop}

\begin{proof}[Proof for R-branes]
	We use Lemma \ref{b min curve}, and parts (1) and (2) of Lemma \ref{sign of b min curve}. Both compositions in the proposition are isotopic to composing a map with a continuation map. In the first composition, there is no contribution from \(\epsilon(\gamma)\) for any loops \(\gamma\) on \(X,Y,Z\) since the representations send all loops on \(X,Y,Z\) to \(1\); this is automatic if \(\dim(M)\geq2\), and in the \(1\)-dimensional case this holds since we imposed that the spin structure on the circle is null-cobordant and so the generator of \(\pi_1(S^1)\) gets sent to \(1\) by the representation (see \cites{Leverson,CN}). In the second composition, there is the contribution of the path from the barycentre of \(s\) to the barycentre of \(t\), which must be accounted for.
\end{proof}

\begin{proof}[Proof for S-branes]
	Similar to \(R\)-branes, both compositions in the proposition are isotopic to composing a map with a continuation map. In the first composition, by looking at the intersection region of the supports of the 3 sheaves involved, we see that the composition is simply multiplication. The second composition is exactly the structure maps of the local system on \(\Sigma\), viewed as a constructible sheaf on the Whitney triangulation \(K\).
\end{proof}

\subsection{Proof of the main theorem}

In this section, we prove:

\begin{thm}\label{main theorem}
	Let \(M\) be a real analytic manifold. Let \(\mathfrak{f}\) be a subanalytic singular isotropic stop bounded in the vertical direction. If \(M\) is 4- or 5-dimensional, we additionally assume Assumption \ref{assumption}. Then, there is a Morita equivalence \[\mathcal{V}(J^1M,\mathfrak{f})\cong \Sh(J^1M,\mathfrak{f}).\]
  	
\end{thm}

\begin{prop}\label{enlarged triangulated equivalence}
	There is a Morita equivalence \[\mathcal{V}(J^1M,\mathfrak{S}\cup\mathfrak{g})/\mathcal{P}\cong\Sh(J^1M,\mathfrak{S}\cup\mathfrak{g})/\mathcal{P}.\]

	Additionally, under this equivalence, \(\Sigma\) as an \(R\)-brane is identified with \(\Sigma\) as an \(S\)-brane, and a linking sphere as an \(R\)-brane is identified with a linking sphere as an \(S\)-brane.
\end{prop}

\begin{proof}
   By Proposition \ref{quotient category generators}, both categories have a set of generators \(\{T_s:s\in\mathcal{S}\}\cup\{Q_\Delta:\Delta\in K\}\). By Proposition \ref{morphisms between generators}, between generators, morphisms are concentrated in degree \(0\). Therefore, after restricting to the full subcategory of the generators, both categories are formal (quasi-equivalent to their cohomology categories). By Proposition \ref{compositions between generators}, after restricting to the full subcategories of the generators, the two categories are equivalent. This proves the first part of the proposition.
	
	By the Yoneda lemma, the fact that the \(R\)-brane \(\Sigma\) is identified with the \(S\)-brane \(\Sigma\) is a consequence of the second part of Proposition \ref{compositions between generators}.

	A similar Yoneda lemma argument with morphism computations show that a linking sphere as an \(R\)-brane is identified with a linking sphere as an \(S\)-brane (linking spheres are not in bijection with strata, but all the arguments of Proposition \ref{morphisms between generators} and Proposition \ref{compositions between generators} still apply, since morphisms from linking spheres to generators only depend on the isotopy type of a linking sphere union a generator after sufficiently wrapping the generator, which turns out to always be isotopic to an unknot and a pushoff of it).
\end{proof}

Finally, we prove the main theorem:

\begin{proof}[Proof of Theorem \ref{main theorem}]
	Let \(\mathcal{C}\) denote either \(\mathcal{V}\) or \(\Sh\). We have that \(\mathcal{C}(J^1M,\mathfrak{f})\) is a quotient category of a subcategory of \(\mathcal{C}(J^1M,\mathfrak{S}\cup\mathfrak{g})/\mathcal{P}\). The subcategory of \(\mathcal{C}(J^1M,\mathfrak{S}\cup\mathfrak{g})/\mathcal{P}\) is equivalent to \(\mathcal{C}(J^1M,\mathfrak{S})\), generated by \(\{T_{s}:s\in \mathcal{S}\}\cup\{\Sigma\text{ below }T^*M\times \{0\}\}\). Since the equivalence in Proposition \ref{enlarged triangulated equivalence} sends \(\Sigma\) as an \(R\)-brane to \(\Sigma\) as an \(S\)-brane, the two subcategories \(\mathcal{C}(J^1M,\mathfrak{S})\) are Morita equivalent. To obtain \(\mathcal{C}(J^1M,\mathfrak{f})\), we need to take the quotient of \(\mathcal{C}(J^1M,\mathfrak{S})\) by all linking spheres of \(\mathfrak{S}\setminus \mathfrak{f}\). Since the equivalence in proposition \ref{enlarged triangulated equivalence} sends linking spheres as \(R\)-branes to linking spheres as \(S\)-branes, we have that the two categories \(\mathcal{C}(J^1M,\mathfrak{f})\) are Morita equivalent, which finishes the proof.
\end{proof}

\appendix

\section{Some holomorphic curves and their signs}\label{sign appendix}

In this appendix, we use the technique of quantum flow trees, developed in \cite{EESa}, \cite{EENS}, \cite{EL}. We refer the reader to \cite{Liu25}*{Section 2.2} for an exposition that suffices for the applications we need here.

Let \(\Lambda\subset J^1M\) be a Legendrian \(R\)-brane. We identify a contact neighborhood of \(\Lambda\) with a contact neighborhood \(U\) of the \(0\)-section in \(J^1\Lambda\). Let \(f,h:\Lambda\to \R\) be small positive functions. Let
\[
f_1=0, f_2=f, f_3=f+\epsilon f+\epsilon^2 h,
\]
Let \(\Lambda_i:=\Gamma_{df_i,\Lambda}\subset U\) be the \(1\)-jet graph of \(f_i\). We would like to degenerate \(\Lambda_1\cup\Lambda_2\cup\Lambda_3\) onto \(\Lambda_1\cup\Lambda_3\). Perturb \(f,h\) so that \(\Lambda_1\cup\Lambda_2\cup\Lambda_3\) is quantum flow tree-generic. Equip \(\Lambda_i\) with the structure of a Legendrian \(R\)-brane coming from that of \(\Lambda\) (the brane structure won't be present until Proposition \ref{c is a cocycle}).

We assume that the function \(f\) involved in the perturbation scheme admits only 1 critical point of index 0. We write \(\min_{ij}\) for the short Reeb chord corresponding to the critical points of index 0 from \(\Lambda_i\) to \(\Lambda_j\). By \cite{Liu25}*{Corollary 2.12}, we see that:

\begin{lem}\label{min min curve}
	There is a unique holomorphic curve \(\Delta_m\) with positive puncture at \(\min_{13}\), and with negative punctures at \(\min_{23},\min_{12}\).
\end{lem}

Since each \(R\)-brane structures on all \(\Lambda_i\) correspond to each other, the orientation of the holomorphic curve described in Lemma \ref{min min curve} depends only on \(n=\dim(\Lambda)\). We write \((-1)^\dagger\) for the orientation of this curve, where \(\dagger\) is a \(\Z/2\)-valued function of \(n\).

\begin{rem}
   In this paper, we do not explicitly compute the value of \(\dagger\). It seems likely that with appropriate choice, \(\dagger=0\) and does not depend on \(n\), but the author does not know how to prove this.
\end{rem}

Our next goal is to determine the signs of the following curves: Let \(\Lambda_0\) be a Legendrian \(R\)-brane which is disjoint from \(\Lambda_1,\Lambda_2,\Lambda_3\). Isotope \(\Lambda_1,\Lambda_2,\Lambda_3\) so that: they are close enough, and for \(i=1,2\), \(\Lambda_0\cup\Lambda_i\cup\Lambda_{i+1}\) is quantum flow tree generic for the purpose of degenrating onto \(\Lambda_0\cup\Lambda_i\). Let \(1\leq i,j\leq3\). We have that the set of Reeb chords \(R_{\Lambda_0,\Lambda_i}\) (which means Reeb chords from \(\Lambda_0\) and \(\Lambda_i\)) is in bijection with \(R_{\Lambda_0,\Lambda_j}\). Similarly, \(R_{\Lambda_i,\Lambda_0}\) is in bijection with \(R_{\Lambda_j,\Lambda_0}\). Under these bijections, \(b_{0i}\) (resp. \(b_{i0}\)) is identified with \(b_{0j}\) (resp. \(b_{j0}\)). By the exact same argument as in \cite{Liu25}*{Corollary 2.12}, we have:

\begin{lem}\label{b min curve}
	Let \(1\leq i<j\leq 3\). For each chord \(b_{0i}\), there is a unique rigid holomorphic curve \(\Delta_{0ij}\) with positive puncture at \(b_{0,j}\), and negative punctures at \(\min_{i,j}\) and \(b_{0i}\). For each chord \(b_{i0}\), there is a unique rigid holomorphic curve \(\Delta_{ij0}\) with positive puncture at \(b_{i0}\), and negative punctures at \(b_{j,0}\) and \(\min_{i,j}\).
\end{lem}

\begin{rem}
	Since it suffices to state Lemma \ref{b min curve} and Lemma \ref{sign of b min curve} with only 2 copies \(\Lambda_1,\Lambda_2\), it may seem weird to state them using 3 parallel copies of \(\Lambda\). The reason we use 3 copies is that they are necessary for the proof of Lemma \ref{sign of b min curve}.
\end{rem}

\begin{lem}\label{sign of b min curve}
	The signs \(\sigma\) of the curves in Lemma \ref{b min curve} are determined as follows:
	\begin{enumerate}
	\item If \(|b_{0i}|\) is odd, then \(\sigma(\Delta_{0ij})=(-1)^\dagger\).
	\item If \(|b_{i0}|\) is odd, then \(\sigma(\Delta_{ij0})=(-1)^\dagger\).
	\item If \(|b_{0i}|\) is even, then \(\sigma(\Delta_{0ij})=(-1)^{\dagger+n}\).
	\item If \(|b_{i0}|\) is even, then \(\sigma(\Delta_{ij0})=(-1)^{\dagger+n+1}\).
	\end{enumerate}
\end{lem}

\begin{proof}
	Since \(\Lambda_1,\Lambda_2,\Lambda_3\) are close enough to each other, and they are equipped with the same orientation, relative pin structure, and capping paths, we have that the signs of the curves don't depend on the value of \(i\) or \(j\).

	Note that when we are considering the ends of a 1-dimensional moduli space of holomorphic curves, the orientations of the ends are not given by simply multiplying the signs of the rigid curves in the two story holomorphic curve buildings. Instead, the signs of the ends are given by the formula in \cite{EES2}*{Lemma 4.11}.

	There is a 1-dimensional moduli space of holomorphic curves, whose 2 end points are given by two story holomorphic curve buildings as follows: on one end, we glue the positive puncture of \(\Delta_m\) (which is \(\min_{13}\)) to the first negative puncture of \(\Delta_{013}\); on the opposite end, we glue the positive puncture of \(\Delta_{012}\) (which is \(b_{02}\)) to the second negative puncture of \(\Delta_{023}\).
	
	If \(|b_{0i}|\) is odd, then \cite{EES2}*{Lemma 4.11} tells us that 
	\[
	\sigma(\Delta_m)\sigma(\Delta_{013})=\sigma(\Delta_{012})\sigma(\Delta_{023}).
	\]
	Since \(\sigma(\Delta_{013})=\sigma(\Delta_{012})=\sigma(\Delta_{023})\), we have that \(\sigma(\Delta_{013})=\sigma(\Delta_{012})=\sigma(\Delta_{023})=(-1)^\dagger\). This proves (1).

	If \(|b_{0i}|\) is even, then \cite{EES2}*{Lemma 4.11} tells us that
	\[
	\sigma(\Delta_m)\sigma(\Delta_{013})=(-1)^n\sigma(\Delta_{012})\sigma(\Delta_{023}),
	\]
	where \(n\) is the dimension of Legendrians here. Since \(\sigma(\Delta_{013})=\sigma(\Delta_{012})=\sigma(\Delta_{023})\), we have that \(\sigma(\Delta_{013})=\sigma(\Delta_{012})=\sigma(\Delta_{023})=(-1)^{\dagger+n}\). This proves (3).

	There is a 1-dimensional moduli space of holomorphic curves, whose 2 end points are given by two story holomorphic curve buildings as follows: on one end, we glue the positive puncture of \(\Delta_m\) (which is \(\min_{13}\)) to the second negative puncture of \(\Delta_{130}\); on the opposite end, we glue the positive puncture of \(\Delta_{230}\) (which is \(b_{20}\)) to the first negative puncture of \(\Delta_{120}\).

	If \(|b_{i0}|\) is odd, then \cite{EES2}*{Lemma 4.11} tells us that
	\[
	\sigma(\Delta_m)\sigma(\Delta_{130})=\sigma(\Delta_{230})\sigma(\Delta_{120}),
	\]
	which implies \(\sigma(\Delta_{130})=\sigma(\Delta_{230})=\sigma(\Delta_{120})=(-1)^\dagger\). This proves (2).

	If \(|b_{i0}|\) is even, then \cite{EES2}*{Lemma 4.11} tells us that
	\[
	\sigma(\Delta_m)\sigma(\Delta_{130})=(-1)^{n+1}\sigma(\Delta_{230})\sigma(\Delta_{120}),
	\]
	which implies \(\sigma(\Delta_{130})=\sigma(\Delta_{230})=\sigma(\Delta_{120})=(-1)^{\dagger+n+1}\). This proves (4).
\end{proof}

Combining the above with the sign \((-1)^{\sigma+(n-1)(|q^\vee|)}\) in the \(A_\infty\)-operation, and noting that the factors from loops cancel out since the \(\Lambda_i\) have the same representation, we obtain:

\begin{prop}\label{min m_2 calculation}
	The curves \(\Delta_{120},\Delta_{012}\) contribute to \(m_2\) with the following signs:
	\begin{enumerate}
	\item \(m_2(Bb_{01},A\min_{12})=(-1)^\dagger BAb_{02}\),
	\item \(m_2(A\min_{12},Bb_{20})=(-1)^{\dagger+|b_{10}^\vee|}ABb_{10}\).
	\end{enumerate}
\end{prop}

Next we turn our attention back to curves on \(\Lambda_1\cup\Lambda_2\cup\Lambda_3\). For each Reeb chord \(q\) on \(\Lambda\), there are corresponding long Reeb chords \(q_{ij}\) for \(1\leq i,j\leq 3\). By \cite{Liu25}*{Lemma 2.11, Corollary 2.12}, we have the following:

\begin{lem}
	There is a unique rigid holomorphic curve of the following form:
	\begin{enumerate}
	\item Positive puncture at \(q_{13}\), negative punctures at \(\min_{23},q_{12}\), denoted \(\Delta_1\).
	\item Positive puncture at \(q_{13}\), negative punctures at \(q_{23},\min_{12}\), denoted \(\Delta_2\).
	\item Positive puncture at \(q_{12}\), negative punctures at \(q_{22},\min_{12}\), denoted \(\Delta_3\).
	\item Positive puncture at \(q_{23}\), negative punctures at \(\min_{23},q_{33}\), denoted \(\Delta_4\).
	\end{enumerate}
\end{lem}

By Lemma \ref{sign of b min curve}, we have that if \(|q|\) is even, then \(\sigma(\Delta_1)=(-1)^{\dagger+n}\), and \(\sigma(\Delta_2)=(-1)^{\dagger+n+1}\). We can thus deduce:

\begin{lem}\label{lemma quantum flow trees cancel}
	If \(|q|\) is even, then \(\sigma(\Delta_3)=-\sigma(\Delta_4)\).
\end{lem}
\begin{proof}
	There is a 1-dimensional moduli space of holomorphic curves, whose 2 end points are given by two story holomorphic curve buildings as follows: on one end, we glue the positive puncture of \(\Delta_3\) (which is \(q_{12}\)) to the second negative puncture of \(\Delta_{1}\); on the opposite end, we glue the positive puncture of \(\Delta_4\) (which is \(q_{23}\)) to the first negative puncture of \(\Delta_{2}\). Then, \cite{EES2}*{Lemma 4.11} tells us that
	\[
	\sigma(\Delta_1)\sigma(\Delta_3)=\sigma(\Delta_2)\sigma(\Delta_4),
	\]
	which concludes the proof, since we know \(\sigma(\Delta_1)=-\sigma(\Delta_2)\).
\end{proof}

\begin{rem}
	Similarly, if \(|q|\) is odd, then \(\sigma(\Delta_3)=\sigma(\Delta_4)\), but we do not need this fact.
\end{rem}

Let \(c\) be an index 1 critical point of \(f\). For \(1\leq i<j\leq3\), let \(c_{ij}\) denote the corresponding short Reeb chord from \(\Lambda_i\) to \(\Lambda_j\). By \cite{Liu25}*{Corollary 2.12}, there is a unique rigid holomorphic curve \(\Delta_{c}\) with positive puncture at \(c_{13}\), and negative punctures at \(\min_{23},c_{12}\); and there is a unique rigid holomorphic curve \(\Delta_c '\) with positive puncture at \(c_{13}\), and negative punctures at \(c_{23},\min_{12}\). There are exactly 2 rigid holomorphic curves \(L_{ij},L_{ij}'\) with positive puncture at \(c_{ij}\) and negative puncture at \(\min_{ij}\), corresponding to the 2 Morse flow lines.

\begin{lem}\label{lemma flow lines cancel}
	\(\sigma(L_{ij})=-\sigma(L_{ij}')\).
\end{lem}
\begin{proof}
	Similar to before, \(\sigma(L_{ij})\) and \(\sigma(L_{ij}')\) don't depend on the values of \(i,j\).

	There is a 1-dimensional moduli space of holomorphic curves, whose 2 end points are given by two story holomorphic curve buildings as follows: on one end, we glue the positive puncture of \(L_{12}\) (which is \(c_{12}\)) to the second negative puncture of \(\Delta_c\); on the opposite end, we glue the positive puncture of \(\Delta_m\) (which is \(\min_{13}\)) to the negative puncture of either \(L_{13}\) (it is one of \(L_{13}\) or \(L_{13}'\), but we can choose the labeling of these \(L\)'s so that the one here is \(L_{13}\)). Then, \cite{EES2}*{Lemma 4.11} tells us that
	\[
	\sigma(L_{12})\sigma(\Delta_c)=(-1)^{n+1}\sigma(\Delta_m)\sigma(L_{13}).
	\]
	Similarly,
	\[
	\sigma(L_{23}')\sigma(\Delta_{c}')=(-1)^{n}\sigma(\Delta_m)\sigma(L_{13}).
	\]
	So, \(\sigma(L_{12})=-\sigma(L_{23}')\), which concludes the proof.
\end{proof}

\begin{prop}\label{c is a cocycle}
	\(I\min_{12}^\vee\in LCC^0(\Lambda_1,\Lambda_2)\) is a cocycle, that is, \(m_1(I\min_{12}^\vee)=0\).
\end{prop}
\begin{proof}
	By \cite{Liu25}*{Lemma 2.11}, the rigid holomorphic curves on \(\Lambda_1\cup\Lambda_2\) with a negative puncture at \(\min_{12}\) are of the form \(L_{12},L_{12}',\Delta_3,\Delta_4\) (for each choice of an index 1 critical point of \(f\), there is corresponding \(L_{12},L_{12}'\), and for each choice of a Reeb chord \(q\) on \(\Lambda\), there is corresponding \(\Delta_3,\Delta_4\)).
   
   For any index 1 critical point \(c\) of \(f_{i+1}-f_{i}\), let \(\gamma_c\) denote the loop on \(\Pi(\Lambda)\) formed by the capping path from the basepoint to \(c\), the flow line from \(c\) to \(\min\), and the inverse of the capping path from the basepoint to \(\min\). We have that the curves \(L_{12}\) and \(L'_{12}\) each contribute \(\pm\epsilon^\Lambda(\gamma_c)I\epsilon^\Lambda(\gamma_c^{-1})q^\vee=\pm I c^\vee\). Then, Lemma \ref{lemma flow lines cancel} tells us that the summands of \(m_1(I\min_{12}^\vee)\) coming from \(L_{12},L'_{12}\) (for all \(c\)) add up to \(0\).
   
   Let \(\gamma_b\) denote the loop on \(\Pi(\Lambda)\) formed by the capping path from the basepoint to \(\min\), the inverse of the partial flow line from the positive endpoint of \(b\) to \(\min\), and the inverse of the capping path from the basepoint to the positive endpoint of \(b\). Also, let \(\gamma_b'\) denote the loop on \(\Lambda\) formed by the capping path from the basepoint to the negative endpoint of \(b\), the partial flow line from the negative endpoint of \(b\) to \(\min\), and the inverse of the capping path from the basepoint to \(\min\). Then, in \(m_1(I\min_{12}^\vee)\), \(\Delta_4\) contributes the term \(\pm\epsilon(b)\epsilon(\gamma_b)\epsilon(\gamma_b^{-1})b_{i,i+1}^{\vee}=\pm\epsilon(b)b_{i,i+1}^{\vee}\), and \(\Delta_3\) contributes the summand \(\pm\epsilon(b)\epsilon(\gamma_b')\epsilon({\gamma_b'}^{-1})b_{i,i+1}^{\vee}=\pm\epsilon(b)b_{i,i+1}^{\vee}\). Lemma \ref{lemma quantum flow trees cancel} then tells us that the summands of \(m_1(I\min_{12}^\vee)\) coming from \(\Delta_3,\Delta_4\) (for all \(b\)) add up to \(0\).
\end{proof}

\begin{proof}[Proof of Proposition \ref{unitality of c}]
	All statements in this appendix generalize in a straightforward way to the case of \(f\) having an arbitrary number of index 0 critical points. Alternatively we can simply define \(c\) by first isotoping \(\Lambda_1\cup\Lambda_2\) so that \(f\) only has 1 critical point of index 0.

	The algebra structure on \(LCH(\Lambda)\) is given by the cohomology of \(m_2\): 
	\[
	H^*m_2:LCH(\Lambda_1,\Lambda_2)\otimes LCH(\Lambda_2,\Lambda_3)\to LCH(\Lambda_1,\Lambda_3)
	\]
	and identifying all \(LCH(\Lambda_i,\Lambda_j)\) with \(LCH(\Lambda)\). To see that \(H^*m_2\) is a well-defined bilinear multiplication map, see the discussion following \cite{Liu25}*{Proposition 3.12}. To see that the multiplication is associative, we need to add one more parallel copy \(\Lambda_4\), and use the fact that \(A_\infty\)-multiplication is associative after taking cohomology.
	
	To show that \(c^\vee\) is a unit, we need to first show that it is an element of \(LCH(\Lambda)\), that is, it is a cocycle. This is Proposition \ref{c is a cocycle}.
	
	We then need to show that \(c^\vee\) acts as the identity in multiplication. By the sign convention in \cite{Seidel}, when taking cohomology, the multiplication is given by \([a]\cdot[b]=(-1)^{|b|}[m_2(a,b)]\), we have by Proposition \ref{min m_2 calculation} that \(c^\vee\) acts as the unit, both in the multiplication on \(LCH(\Lambda)\), and on the modules \(LCH(\Lambda,\Lambda'),LCH(\Lambda',\Lambda)\).
\end{proof}

\begin{lem}\label{isotopy sends c to c}
	Given an isotopy from \(\Lambda\) to \(\Lambda'\), the induced map \[LCH(\Lambda)\to LCH(\Lambda')\] sends \(c^\vee\) to \((c')^\vee\).
\end{lem}
\begin{proof}
	This is because an algebra map sends the unit to the unit.
\end{proof}

\subsection*{Technology disclosure}
In order to draw vector graphics in this paper, the author asked ChatGPT how to download Inkscape when the Inkscape website was down.
\addcontentsline{toc}{section}{Technology disclosure}

\addcontentsline{toc}{section}{References}
\begin{bibdiv}
  \begin{biblist}

\bib{Abouzaid}{article}{
   author={Abouzaid, Mohammed},
   title={A geometric criterion for generating the Fukaya category},
   journal={Publ. Math. Inst. Hautes \'Etudes Sci.},
   number={112},
   date={2010},
   pages={191--240},
   issn={0073-8301},
   review={\MR{2737980}},
   doi={10.1007/s10240-010-0028-5},
}

\bib{Abouzaid2}{article}{
   author={Abouzaid, Mohammed},
   title={A cotangent fibre generates the Fukaya category},
   journal={Adv. Math.},
   volume={228},
   date={2011},
   number={2},
   pages={894--939},
   issn={0001-8708},
   review={\MR{2822213}},
   doi={10.1016/j.aim.2011.06.007},
}

\bib{Abouzaid3}{article}{
   author={Abouzaid, Mohammed},
   title={Nearby Lagrangians with vanishing Maslov class are homotopy
   equivalent},
   journal={Invent. Math.},
   volume={189},
   date={2012},
   number={2},
   pages={251--313},
   issn={0020-9910},
   review={\MR{2947545}},
   doi={10.1007/s00222-011-0365-0},
}

\bib{Abouzaid4}{article}{
   author={Abouzaid, Mohammed},
   title={On the wrapped Fukaya category and based loops},
   journal={J. Symplectic Geom.},
   volume={10},
   date={2012},
   number={1},
   pages={27--79},
   issn={1527-5256},
   review={\MR{2904032}},
   doi={10.4310/jsg.2012.v10.n1.a3},
}

\bib{ABS}{article}{
   author={An, Byung Hee},
   author={Bae, Youngjin},
   author={Su, Tao},
   title={Augmentations are sheaves for Legendrian graphs},
   journal={J. Symplectic Geom.},
   volume={20},
   date={2022},
   number={2},
   pages={259--415},
   issn={1527-5256},
   review={\MR{4537126}},
   doi={10.4310/jsg.2022.v20.n2.a1},
}

\bib{BC}{article}{
   author={Bourgeois, Fr\'ed\'eric},
   author={Chantraine, Baptiste},
   title={Bilinearized Legendrian contact homology and the augmentation
   category},
   journal={J. Symplectic Geom.},
   volume={12},
   date={2014},
   number={3},
   pages={553--583},
   issn={1527-5256},
   review={\MR{3248668}},
   doi={10.4310/jsg.2014.v12.n3.a5},
}

\bib{CG1}{article}{
   author={Casals, Roger},
   author={Gao, Honghao},
   title={Infinitely many Lagrangian fillings},
   journal={Ann. of Math. (2)},
   volume={195},
   date={2022},
   number={1},
   pages={207--249},
   issn={0003-486X},
   review={\MR{4358415}},
   doi={10.4007/annals.2022.195.1.3},
}

\bib{CG2}{article}{
   author={Casals, Roger},
   author={Gao, Honghao},
   title={A Lagrangian filling for every cluster seed},
   journal={Invent. Math.},
   volume={237},
   date={2024},
   number={2},
   pages={809--868},
   issn={0020-9910},
   review={\MR{4768635}},
   doi={10.1007/s00222-024-01268-y},
}

\bib{CN}{article}{
   author={Casals, Roger},
   author={Ng, Lenhard},
   title={Braid loops with infinite monodromy on the Legendrian contact DGA},
   journal={J. Topol.},
   volume={15},
   date={2022},
   number={4},
   pages={1927--2016},
   issn={1753-8416},
   review={\MR{4584583}},
   doi={10.1112/topo.12264},
}

\bib{CDGG}{article}{
   author={Chantraine, Baptiste},
   author={Dimitroglou Rizell, Georgios},
   author={Ghiggini, Paolo},
   author={Golovko, Roman},
   title={Noncommutative augmentation categories},
   conference={
      title={Proceedings of the G\"okova Geometry-Topology Conference 2015},
   },
   book={
      publisher={G\"okova Geometry/Topology Conference (GGT), G\"okova},
   },
   isbn={978-1-57146-321-0},
   date={2016},
   pages={116--150},
   review={\MR{3526840}},
}

\bib{CDGG2}{article}{
   author={Chantraine, Baptiste},
   author={Dimitroglou Rizell, Georgios},
   author={Ghiggini, Paolo},
   author={Golovko, Roman},
   title={Geometric generation of the wrapped Fukaya category of Weinstein
   manifolds and sectors},
   journal={Ann. Sci. \'Ec. Norm. Sup\'er. (4)},
   volume={57},
   date={2024},
   number={1},
   pages={1--85},
   issn={0012-9593},
   review={\MR{4732675}},
   doi={10.24033/asens.2570},
}



\bib{CNS}{article}{
   author={Chantraine, Baptiste},
   author={Ng, Lenhard},
   author={Sivek, Steven},
   title={Representations, sheaves and Legendrian $(2,m)$ torus links},
   journal={J. Lond. Math. Soc. (2)},
   volume={100},
   date={2019},
   number={1},
   pages={41--82},
   issn={0024-6107},
   review={\MR{3999682}},
   doi={10.1112/jlms.12204},
}

\bib{Che}{article}{
   author={Chekanov, Yuri},
   title={Differential algebra of Legendrian links},
   journal={Invent. Math.},
   volume={150},
   date={2002},
   number={3},
   pages={441--483},
   issn={0020-9910},
   review={\MR{1946550}},
   doi={10.1007/s002220200212},
}

\bib{Cza}{article}{
   author={Czapla, Ma\l gorzata},
   title={Definable triangulations with regularity conditions},
   journal={Geom. Topol.},
   volume={16},
   date={2012},
   number={4},
   pages={2067--2095},
   issn={1465-3060},
   review={\MR{3033513}},
   doi={10.2140/gt.2012.16.2067},
}

\bib{Ekh}{article}{
   author={Ekholm, Tobias},
   title={Morse flow trees and Legendrian contact homology in 1-jet spaces},
   journal={Geom. Topol.},
   volume={11},
   date={2007},
   pages={1083--1224},
   issn={1465-3060},
   review={\MR{2326943}},
   doi={10.2140/gt.2007.11.1083},
}

\bib{EESa}{article}{
   author={Ekholm, Tobias},
   author={Etnyre, John B.},
   author={Sabloff, Joshua M.},
   title={A duality exact sequence for Legendrian contact homology},
   journal={Duke Math. J.},
   volume={150},
   date={2009},
   number={1},
   pages={1--75},
   issn={0012-7094},
   review={\MR{2560107}},
   doi={10.1215/00127094-2009-046},
}

\bib{EES1}{article}{
   author={Ekholm, Tobias},
   author={Etnyre, John},
   author={Sullivan, Michael},
   title={The contact homology of Legendrian submanifolds in ${\Bbb
   R}^{2n+1}$},
   journal={J. Differential Geom.},
   volume={71},
   date={2005},
   number={2},
   pages={177--305},
   issn={0022-040X},
   review={\MR{2197142}},
}

\bib{EES2}{article}{
   author={Ekholm, Tobias},
   author={Etnyre, John},
   author={Sullivan, Michael},
   title={Orientations in Legendrian contact homology and exact Lagrangian
   immersions},
   journal={Internat. J. Math.},
   volume={16},
   date={2005},
   number={5},
   pages={453--532},
   issn={0129-167X},
   review={\MR{2141318}},
   doi={10.1142/S0129167X05002941},
}

\bib{EES3}{article}{
   author={Ekholm, Tobias},
   author={Etnyre, John},
   author={Sullivan, Michael},
   title={Legendrian contact homology in $P\times\Bbb R$},
   journal={Trans. Amer. Math. Soc.},
   volume={359},
   date={2007},
   number={7},
   pages={3301--3335},
   issn={0002-9947},
   review={\MR{2299457}},
   doi={10.1090/S0002-9947-07-04337-1},
}




\bib{EENS}{article}{
   author={Ekholm, Tobias},
   author={Etnyre, John B.},
   author={Ng, Lenhard},
   author={Sullivan, Michael G.},
   title={Knot contact homology},
   journal={Geom. Topol.},
   volume={17},
   date={2013},
   number={2},
   pages={975--1112},
   issn={1465-3060},
   review={\MR{3070519}},
   doi={10.2140/gt.2013.17.975},
}



\bib{EL}{article}{
   author={Ekholm, Tobias},
   author={Lekili, Yank\i},
   title={Duality between Lagrangian and Legendrian invariants},
   journal={Geom. Topol.},
   volume={27},
   date={2023},
   number={6},
   pages={2049--2179},
   issn={1465-3060},
   review={\MR{4634745}},
   doi={10.2140/gt.2023.27.2049},
}


\bib{EGH}{article}{
   author={Eliashberg, Y.},
   author={Givental, A.},
   author={Hofer, H.},
   title={Introduction to symplectic field theory},
   note={GAFA 2000 (Tel Aviv, 1999)},
   journal={Geom. Funct. Anal.},
   date={2000},
   pages={560--673},
   issn={1016-443X},
   review={\MR{1826267}},
   doi={10.1007/978-3-0346-0425-3\_4},
}

\bib{EHS}{article}{
   author={Eliashberg, Y.},
   author={Hofer, H.},
   author={Salamon, D.},
   title={Lagrangian intersections in contact geometry},
   journal={Geom. Funct. Anal.},
   volume={5},
   date={1995},
   number={2},
   pages={244--269},
   issn={1016-443X},
   review={\MR{1334868}},
   doi={10.1007/BF01895668},
}

\bib{Floer1}{article}{
   author={Floer, Andreas},
   title={Morse theory for fixed points of symplectic diffeomorphisms},
   journal={Bull. Amer. Math. Soc. (N.S.)},
   volume={16},
   date={1987},
   number={2},
   pages={279--281},
   issn={0273-0979},
   review={\MR{0876964}},
   doi={10.1090/S0273-0979-1987-15517-0},
}

\bib{Floer2}{article}{
   author={Floer, Andreas},
   title={Morse theory for Lagrangian intersections},
   journal={J. Differential Geom.},
   volume={28},
   date={1988},
   number={3},
   pages={513--547},
   issn={0022-040X},
   review={\MR{0965228}},
}

\bib{Floer3}{article}{
   author={Floer, Andreas},
   title={A relative Morse index for the symplectic action},
   journal={Comm. Pure Appl. Math.},
   volume={41},
   date={1988},
   number={4},
   pages={393--407},
   issn={0010-3640},
   review={\MR{0933228}},
   doi={10.1002/cpa.3160410402},
}

\bib{Floer4}{article}{
   author={Floer, Andreas},
   title={The unregularized gradient flow of the symplectic action},
   journal={Comm. Pure Appl. Math.},
   volume={41},
   date={1988},
   number={6},
   pages={775--813},
   issn={0010-3640},
   review={\MR{0948771}},
   doi={10.1002/cpa.3160410603},
}


\bib{Fukaya}{article}{
   author={Fukaya, Kenji},
   title={Morse homotopy, $A^\infty$-category, and Floer homologies},
   conference={
      title={Proceedings of GARC Workshop on Geometry and Topology '93},
      address={Seoul},
      date={1993},
   },
   book={
      series={Lecture Notes Ser.},
      volume={18},
      publisher={Seoul Nat. Univ., Seoul},
   },
   date={1993},
   pages={1--102},
   review={\MR{1270931}},
}














\bib{GPS1}{article}{
   author={Ganatra, Sheel},
   author={Pardon, John},
   author={Shende, Vivek},
   title={Covariantly functorial wrapped Floer theory on Liouville sectors},
   journal={Publ. Math. Inst. Hautes \'Etudes Sci.},
   volume={131},
   date={2020},
   pages={73--200},
   issn={0073-8301},
   review={\MR{4106794}},
   doi={10.1007/s10240-019-00112-x},
}

\bib{GPS2}{article}{
   author={Ganatra, Sheel},
   author={Pardon, John},
   author={Shende, Vivek},
   title={Sectorial descent for wrapped Fukaya categories},
   journal={J. Amer. Math. Soc.},
   volume={37},
   date={2024},
   number={2},
   pages={499--635},
   issn={0894-0347},
   review={\MR{4695507}},
   doi={10.1090/jams/1035},
}

\bib{GPS3}{article}{
   author={Ganatra, Sheel},
   author={Pardon, John},
   author={Shende, Vivek},
   title={Microlocal Morse theory of wrapped Fukaya categories},
   journal={Ann. of Math. (2)},
   volume={199},
   date={2024},
   number={3},
   pages={943--1042},
   issn={0003-486X},
   review={\MR{4740209}},
   doi={10.4007/annals.2024.199.3.1},
}

\bib{Gao}{article}{
author={Gao, Honghao},
title={Augmentations and sheaves for links, \emph{arXiv preprint: arXiv:2109.01285}},
journal={},
date={2021},
volume={},
pages={}
}

\bib{GSW}{article}{
   author={Gao, Honghao},
   author={Shen, Linhui},
   author={Weng, Daping},
   title={Augmentations, fillings, and clusters},
   journal={Geom. Funct. Anal.},
   volume={34},
   date={2024},
   number={3},
   pages={798--867},
   issn={1016-443X},
   review={\MR{4743511}},
   doi={10.1007/s00039-024-00673-y},
}

\bib{GL}{article}{
author={Gao, Honghao},
author={Liu, Hanming},
title={Isomorphism in the augmentation category, \emph{arXiv preprint: arXiv:2602.12555}},
journal={},
date={2026},
volume={},
pages={}
}


\bib{Gromov}{article}{
   author={Gromov, M.},
   title={Pseudo holomorphic curves in symplectic manifolds},
   journal={Invent. Math.},
   volume={82},
   date={1985},
   number={2},
   pages={307--347},
   issn={0020-9910},
   review={\MR{0809718}},
   doi={10.1007/BF01388806},
}

\bib{GKS}{article}{
   author={Guillermou, St\'ephane},
   author={Kashiwara, Masaki},
   author={Schapira, Pierre},
   title={Sheaf quantization of Hamiltonian isotopies and applications to
   nondisplaceability problems},
   journal={Duke Math. J.},
   volume={161},
   date={2012},
   number={2},
   pages={201--245},
   issn={0012-7094},
   review={\MR{2876930}},
   doi={10.1215/00127094-1507367},
}





\bib{KS}{book}{
   author={Kashiwara, Masaki},
   author={Schapira, Pierre},
   title={Sheaves on manifolds},
   series={Grundlehren der mathematischen Wissenschaften [Fundamental
   Principles of Mathematical Sciences]},
   volume={292},
   note={With a chapter in French by Christian Houzel},
   publisher={Springer-Verlag, Berlin},
   date={1990},
   pages={x+512},
   isbn={3-540-51861-4},
   review={\MR{1074006}},
   doi={10.1007/978-3-662-02661-8},
}

\bib{Kilgore}{article}{
author={Kilgore, Eric},
title={Legendrian non-squeezing via microsheaves},
journal={arXiv:2412.03823},
date={2024},
volume={},
pages={}
}

\bib{Kontsevich}{article}{
   author={Kontsevich, Maxim},
   title={Homological algebra of mirror symmetry},
   conference={
      title={Proceedings of the International Congress of Mathematicians,
      Vol.\ 1, 2},
      address={Z\"urich},
      date={1994},
   },
   book={
      publisher={Birkh\"auser, Basel},
   },
   isbn={3-7643-5153-5},
   date={1995},
   pages={120--139},
   review={\MR{1403918}},
}


\bib{Kuo}{article}{
   author={Kuo, Christopher},
   title={Wrapped sheaves},
   journal={Adv. Math.},
   volume={415},
   date={2023},
   pages={Paper No. 108882, 71},
   issn={0001-8708},
   review={\MR{4543073}},
   doi={10.1016/j.aim.2023.108882},
}

\bib{Leverson}{article}{
   author={Leverson, Caitlin},
   title={Augmentations and rulings of Legendrian knots},
   journal={J. Symplectic Geom.},
   volume={14},
   date={2016},
   number={4},
   pages={1089--1143},
   issn={1527-5256},
   review={\MR{3601885}},
   doi={10.4310/JSG.2016.v14.n4.a5},
}

\bib{Li}{article}{
   author={Li, Wenyuan},
   title={Estimating Reeb chords using microlocal sheaf theory, \emph{arXiv preprint: 2106.04079}},
   date={2021}
}


\bib{Liu25}{article}{
author={Liu, Hanming},
title={Augmentation categories in higher dimensions, \emph{arXiv preprint: 2508.06822}},
date={2025},
volume={},
pages={}
}








\bib{Nadler}{article}{
   author={Nadler, David},
   title={Microlocal branes are constructible sheaves},
   journal={Selecta Math. (N.S.)},
   volume={15},
   date={2009},
   number={4},
   pages={563--619},
   issn={1022-1824},
   review={\MR{2565051}},
   doi={10.1007/s00029-009-0008-0},
}

\bib{NZ}{article}{
   author={Nadler, David},
   author={Zaslow, Eric},
   title={Constructible sheaves and the Fukaya category},
   journal={J. Amer. Math. Soc.},
   volume={22},
   date={2009},
   number={1},
   pages={233--286},
   issn={0894-0347},
   review={\MR{2449059}},
   doi={10.1090/S0894-0347-08-00612-7},
}

\bib{Ng}{article}{
   author={Ng, Lenhard L.},
   title={Computable Legendrian invariants},
   journal={Topology},
   volume={42},
   date={2003},
   number={1},
   pages={55--82},
   issn={0040-9383},
   review={\MR{1928645}},
   doi={10.1016/S0040-9383(02)00010-1},
}



\bib{NRSSZ}{article}{
   author={Ng, Lenhard},
   author={Rutherford, Dan},
   author={Shende, Vivek},
   author={Sivek, Steven},
   author={Zaslow, Eric},
   title={Augmentations are sheaves},
   journal={Geom. Topol.},
   volume={24},
   date={2020},
   number={5},
   pages={2149--2286},
   issn={1465-3060},
   review={\MR{4194293}},
   doi={10.2140/gt.2020.24.2149},
}

\bib{Pardon}{article}{
   author={Pardon, John},
   title={Contact homology and virtual fundamental cycles},
   journal={J. Amer. Math. Soc.},
   volume={32},
   date={2019},
   number={3},
   pages={825--919},
   issn={0894-0347},
   review={\MR{3981989}},
   doi={10.1090/jams/924},
}

\bib{Seidel}{book}{
   author={Seidel, Paul},
   title={Fukaya categories and Picard-Lefschetz theory},
   series={Zurich Lectures in Advanced Mathematics},
   publisher={European Mathematical Society (EMS), Z\"urich},
   date={2008},
   pages={viii+326},
   isbn={978-3-03719-063-0},
   review={\MR{2441780}},
   doi={10.4171/063},
}

\bib{Shende}{article}{
   author={Shende, Vivek},
   title={The conormal torus is a complete knot invariant},
   journal={Forum Math. Pi},
   volume={7},
   date={2019},
   pages={e6, 16},
   review={\MR{4010558}},
   doi={10.1017/fmp.2019.1},
}

\bib{STZ}{article}{
   author={Shende, Vivek},
   author={Treumann, David},
   author={Zaslow, Eric},
   title={Legendrian knots and constructible sheaves},
   journal={Invent. Math.},
   volume={207},
   date={2017},
   number={3},
   pages={1031--1133},
   issn={0020-9910},
   review={\MR{3608288}},
   doi={10.1007/s00222-016-0681-5},
}

\bib{Sylvan}{article}{
   author={Sylvan, Zachary},
   title={On partially wrapped Fukaya categories},
   journal={J. Topol.},
   volume={12},
   date={2019},
   number={2},
   pages={372--441},
   issn={1753-8416},
   review={\MR{3911570}},
   doi={10.1112/topo.12088},
}

\bib{Tamarkin}{article}{
   author={Tamarkin, Dmitry},
   title={Microlocal condition for non-displaceability},
   conference={
      title={Algebraic and analytic microlocal analysis},
   },
   book={
      series={Springer Proc. Math. Stat.},
      volume={269},
      publisher={Springer, Cham},
   },
   isbn={978-3-030-01588-6},
   isbn={978-3-030-01586-2},
   date={2018},
   pages={99--223},
   review={\MR{3903318}},
   doi={10.1007/978-3-030-01588-6\_3},
}




 \end{biblist}
\end{bibdiv}
\end{document}